\documentclass[]{amsart}
\usepackage{color,latexsym,amssymb,amsmath,amsthm}
\usepackage{comment}
\usepackage{graphicx}

\definecolor{orange}{rgb}{1,0.5,0}

\theoremstyle{definition}

\title{Thinking like Archimedes with a 3D printer}
\author{Oliver Knill and Elizabeth Slavkovsky}
\date{Jan 28, 2013}
\address{
        Department of Mathematics \\
        Harvard University \\
        and Harvard Extension school \\
        Cambridge, MA, 02138
        }


\begin{document}
\maketitle

\begin{abstract}
We illustrate Archimedes' method using models produced with 3D printers.
This approach allowed us to create physical proofs of results known to Archimedes 
and to illustrate ideas of a mathematician who is known both for his 
mechanical inventions as well as his breakthroughs in geometry and calculus. 
We use technology from the 21st century to trace intellectual achievements 
from the 3rd century BC. While we celebrate the 2300th birthday of Archimedes 
(287-212 BC) in 2013, we also live in an exciting time, where 3D printing is 
becoming popular and affordable.
\end{abstract}

\section{Introduction}

Archimedes, whose 2300th birthday we celebrate this year,
was a mathematician and inventor who pursued mechanical methods
to develop a pure theory. In admiration for his mathematical discoveries, 
which were often fueled by experiment, we follow his steps by building 
models produced with modern 3D printers.
Archimedes was an early experimental mathematician 
\cite{BBG} who would use practical problems and experiments as a heuristic tool
\cite{Kropp}. He measured the center of mass and made comparisons between the
volumes and surface areas of known and unknown objects by comparing their
integrals, starting to build the initial ideas of calculus.           
He also tackled problems as a practical engineer, building winches,
pumps, and catapults \cite{Geymonat,Heath3,Heath4,Hirshfeld}. \\ 

In mathematics, Archimedes is credited for the invention of the exhaustion
method, a concept that allowed him to compute the area of a disc,
parabola segments, the volume of a sphere, cones and other quadrics,
objects bound by cylinders and planes, like the hoof, or objects obtained
by the intersection of two or three cylinders, or more general solids that we call 
now Archimedean spheres \cite{Pippenger}. Archimedes'
method was an important step towards calculus; his ideas would later
be built upon by the Italian school like Cavalieri and Torricelli, and by
Descartes and Fermat in France. It would eventually be refined by Leibniz and
Newton to become the calculus we know today \cite{Kramer}. 
Archimedes had a different approach to research than Euclid. 
While appreciating Euclid's format in communicating proof, he did not follow
a deductive way but pursued mechanical methods to develop a pure theory.\\

The life of Archimedes also lead to a dramatic story in the history of 
mathematics \cite{Saito}: the event was the discovery of the ``palimpsest", 
a Byzantine manuscript from the 10th century that
had been reused to copy a prayer book three centuries
later. It remained hidden until the mid-nineteenth century \cite{archimedescodex}
when Johannes Heiberg \cite{Heiberg}, a philologist of Greek mathematics, consulted it
in Istanbul. But only in the 21st century, with the help of modern 
computer technology, researchers finally were able to read the 
manuscript completely \cite{Heath2,Netz}. 
This document revealed much about the thinking of Archimedes and gave historians
insights about how he discovered mathematical results. 
It contains ``The Method", of which the palimpsest contains the only known copy.  \\

\begin{figure}
\scalebox{0.25}{\includegraphics{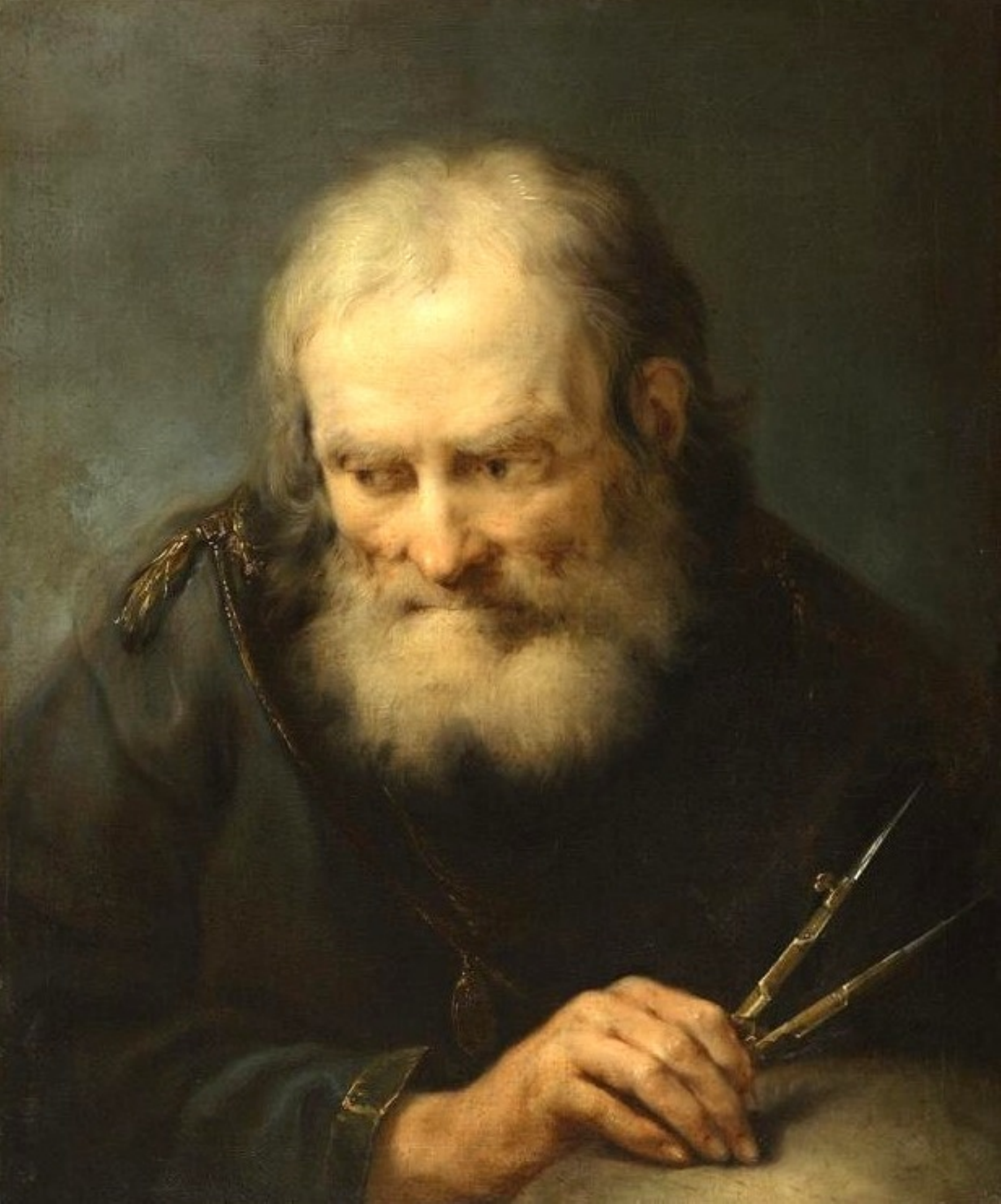}}
\scalebox{0.25}{\includegraphics{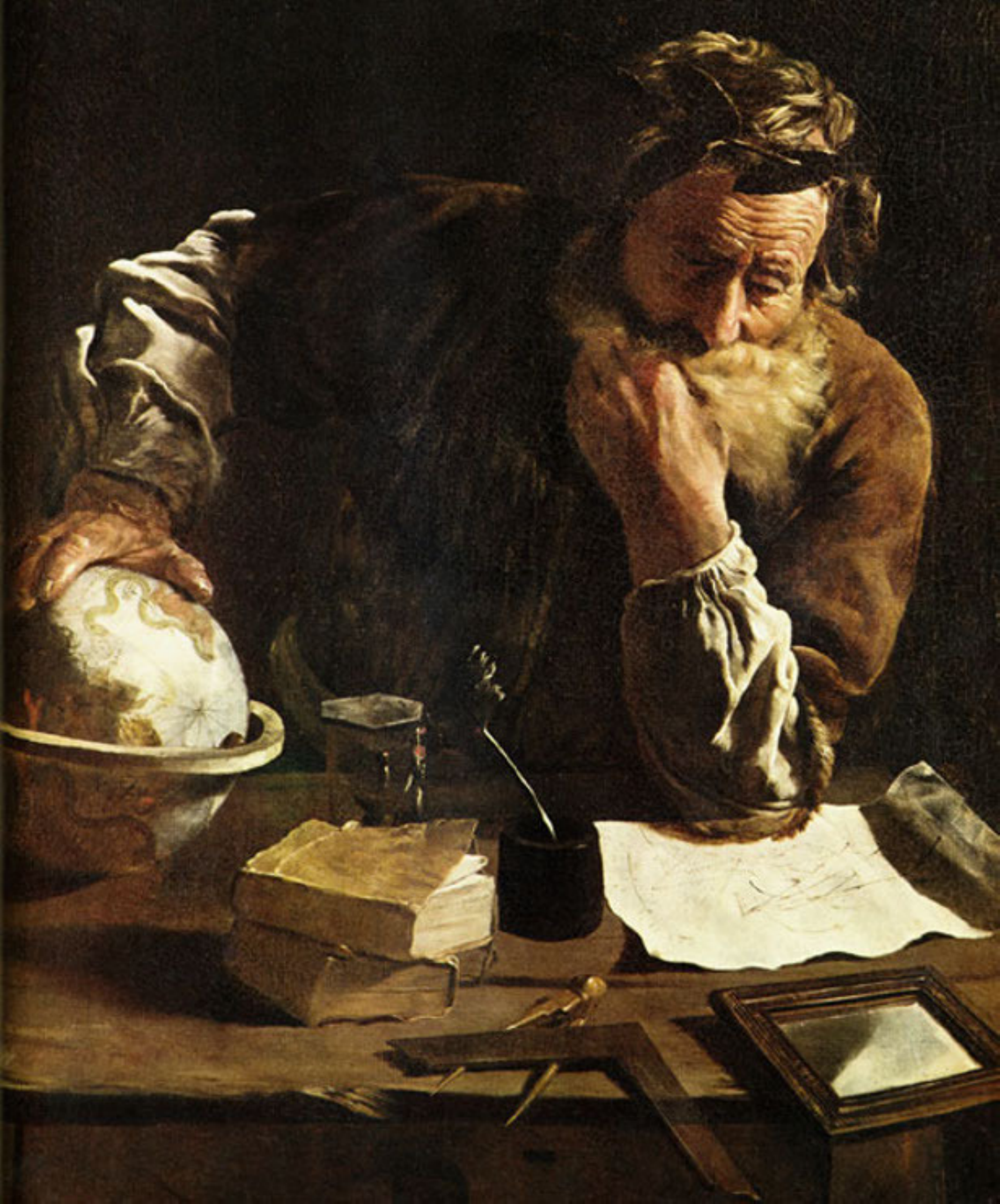}}
\caption{
The left figure shows an oil painting of Archimedes by Giuseppe
Nogari from the mid-18th century.  The oil canvas is located at 
the Pushkin Museum in Moscow. 
The right picture shows a portrait of Archimedes painted around 1620 by
Domenico Fetti in Mantua. It is now located in the art museum
``Alte Meister" in the German town of Dresden \cite{archimedes2}. }
\label{archimedes}
\end{figure}

We know the birth year 287 BC of Archimedes more accurately than from 
other contemporary mathematicians because of a biography about Archimedes 
written by Heraclides \cite{dijksterhuis}. Therefore, the year 2013 celebrates 
the 2300th birthday of Archimedes. 
The death year is known because his city was sacked in 212 BC, an event during 
which Archimedes was killed. We also know much about Archimedes by the work of the historian 
John Tzetzes \cite{Stein}.  

\section{Archimedes mathematics}

Archimedes achievements are covered in all textbooks dealing with the history of 
early geometry. Examples are \cite{Czwalina,Lurje,Turnbull,Becker,Katz2011}. 
While some of the main ideas in calculus were pioneered during Archimedes' time, 
crucial ingredients like the concept of a function or a precise notion of 
limit were missing. Calculus needed to wait 1500 years to be completed \cite{Kramer}. 
Modern revelations have shown that Archimedes already seemed to have made
an important step towards tacking the infinity in calculus \cite{archimedescodex,Netz}.
Only  Newton and Leibniz, after the concept of a ``function" was available, made calculus 
into the science we know today. Despite having modern techniques available to compute,
the method of Archimedes can still be admired today and influenced other techniques in analysis.
His thinking contained the germ for understanding what a limit is. 
The mechanical approach Archimedes used to compute volumes 
consisted of a comparison method: cut a region or body into thin 
slices and then compare the area or volume of slices to different bodies of known area 
or volume. The method allowed him also to compute other integral quantities 
like the center of mass of bodies.   \\

In his work ``On the Sphere and the Cylinder" \cite{Archimedes}, Archimedes rigorously derived 
the surface area and volume of the unit sphere. 
He noticed that the ratio of the surface area $4\pi r^2$ to that of its
circumscribing cylinder is $6 \pi r^2$ is {\bf $2/3$} and that the ratio
the volume of the sphere $4 \pi r^3/3$ with the volume of the circumscribing cylinder
$2\pi r^3$ is {\bf $2/3$} also. He generalized this and also ``Archimedean domes" or 
``Archimedean globes" have volume equal to {\bf $2/3$} of the prism in which they are 
inscribed. It was discovered only later that also for globes, the surface area is {\bf $2/3$}
of the surface area of a circumscribing prism \cite{Apostol}. \\

As reported by Plutarch \cite{Thomas}, Archimedes
was so proud of his sphere computation that he asked to have it inscribed
in his tombstone. Heath concluded from this that Archimedes
regarded it as his greatest achievement. It was Cicero who confirmed that 
``a sphere along with a cylinder had been set up on top of his grave"
\cite{Simms,Jaeger}. Unfortunately, the tombstone is lost. \\

To cite \cite{Dunham}:
{\it ``But for all of these accomplishments, his undisputed masterpiece
was an extensive, two-volume work titled `On the Sphere and the Cylinder'.
Here with almost superhuman cleverness, he determined volumes
and surface areas of spheres and related bodies, thereby achieving for
three-dimensional solids what Measurement of a Circle had done for
two-dimensional figures. It was a stunning triumph, one that Archimedes
himself seems to have regarded as the apex of his career."}
The work ``On the Sphere and the Cylinder" is considered an ideal 
continuation of Euclid's elements \cite{Geymonat,Eves}. \\

The work of Archimedes is now available in \cite{Heath} and in annotated form 
in \cite{Heath2,Netz}. Still, Archimedes can be difficult to read for a modern mathematicians.
Lurje \cite{Lurje} writes on page 176: 
{\it ``Archimedes is a very difficult author. He appears as such for us and must have so for 
ancient mathematicians. If Plutarch lauds the clarity of the proofs of Archimedes, 
then this only shows that
Plutarch did not understand mathematics, that he never read Archimedes and only wanted to paint
a picture of a genius."}
As it happens frequently for creative mathematicians, there could also be
gaps like in Proposition 9 of ``On the Sphere and the Cylinder" \cite{lacuna}. 
But shortcomings are frequent in works of pioneers.
There is no doubt that the mathematics of Archimedes was astounding for the time. 
From \cite{Dunham}: 
{\it ``Archimedes was doing mathematics whose brilliance would be unmatched for
centuries! Not until the development of calculus in the latter years of
the seventeenth century did people advance the understanding of volumes and surface areas
of solids beyond its Archimedean foundation. It is certain that, regardless
of what future glories await the discipline of mathematics, no one
will ever again be 2000 years ahead of his or her time."}

\section{Archimedes engineering}

Also the practical work of Archimedes is impressive.
We dare to speculate that Archimedes - the engineer - would have liked modern 
3D-prototyping technology. It is conceivable that he would have
considered using it as a tool for experimentation or illustration.
Instead of drawing winches, pulleys, or pumps and then build large scale models,
he could prototype smaller versions first to test the device.  \\

Heath writes \cite{Heath}: ``Incidentally he made himself famous
by a variety of ingenious mechanical inventions. These things
were however merely the `diversions of geometry at play' (Plutarch)
and he attached no importance to them." In \cite{Lurje} page 130, this statement is
questioned. In Al Jallil-As-Siisi from the 10th to 11th century for example,
one does not find this ``contemptuousness with respect to practical mechanics".
Lurie claims that Plutarch obviously did not read the work of Archimedes
and cites Tacquet, a geometer from the 17th century: ``Archimedes is more
praised than read, more lauded than understood."
The mechanical inventions include the Cochlias or the Archimedean screw,
which is a machine to pump water, which he may have developed in his youth \cite{Hirshfeld}.
According to Diodorus, it has been used in Egypt for the
irrigation of fields or in Spain in order to pump  water out of mines.
For the description we must rely on Vitruvius \cite{dijksterhuis}.

\begin{figure}
\scalebox{0.20}{\includegraphics{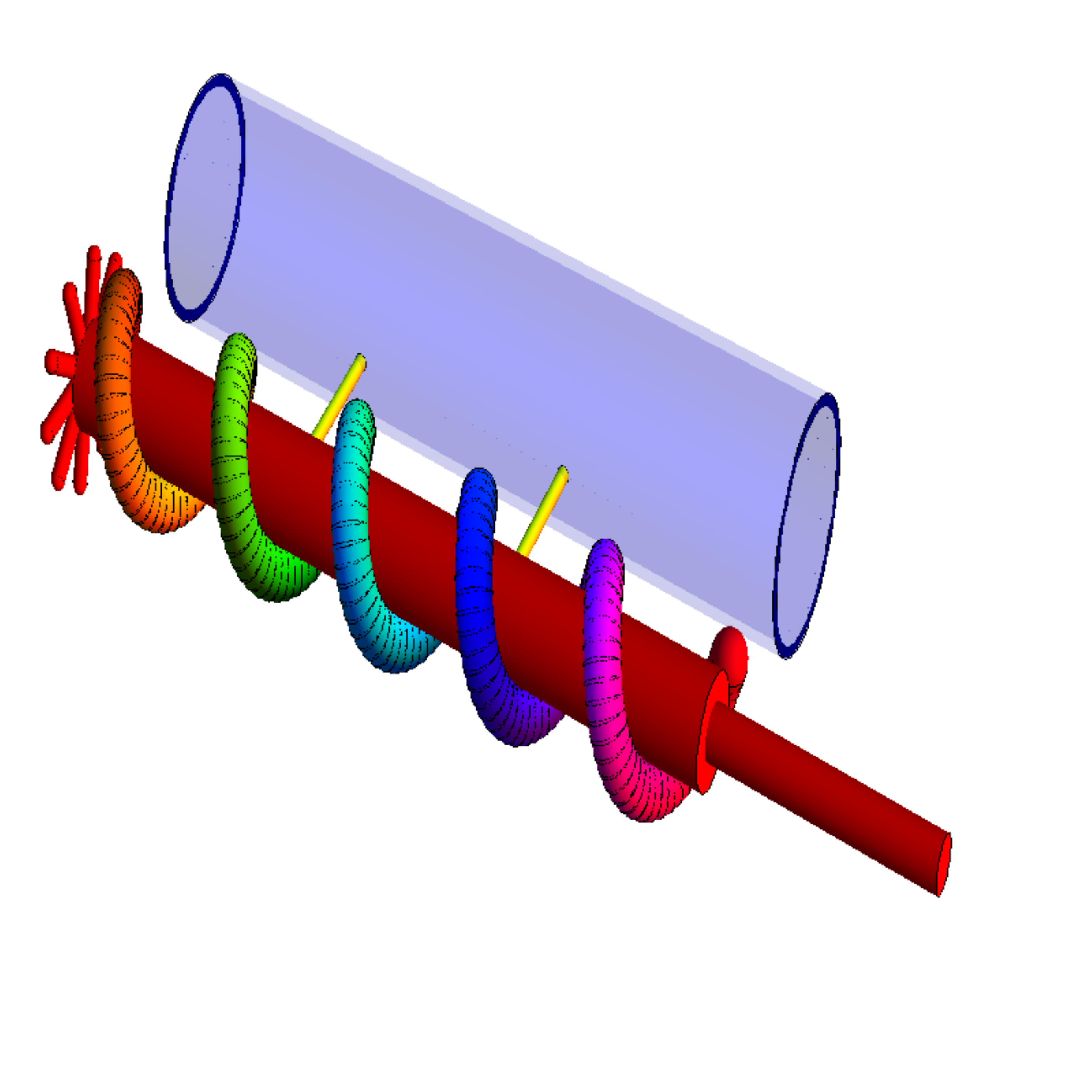}}
\scalebox{0.20}{\includegraphics{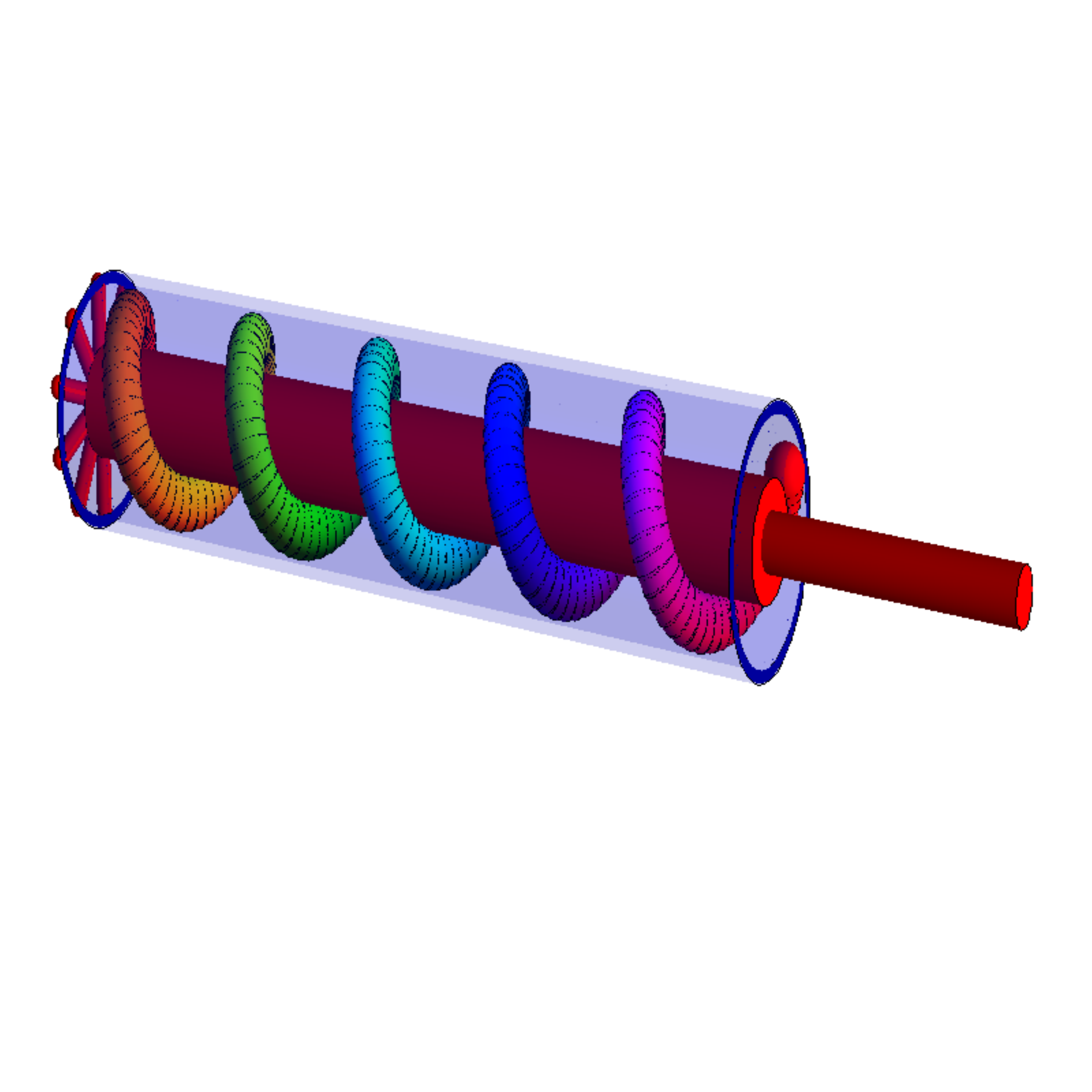}}
\caption{
A model of the ``cochlias" or Archimedes screw. Even so the model is not
a connected body, it is printed with connectors, which are then broken off. 
}
\label{screw}
\end{figure}
\begin{figure}
\scalebox{0.20}{\includegraphics{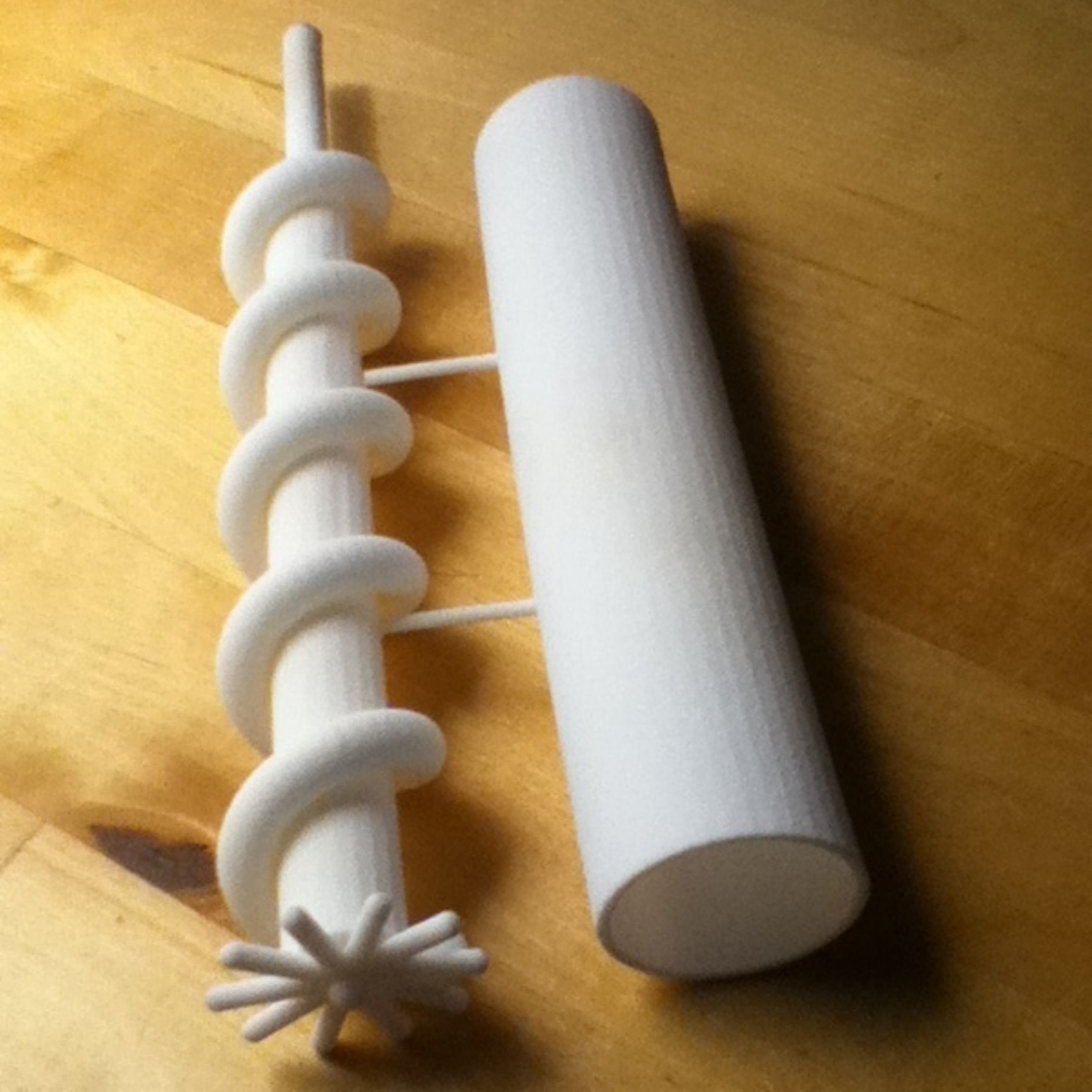}}
\scalebox{0.20}{\includegraphics{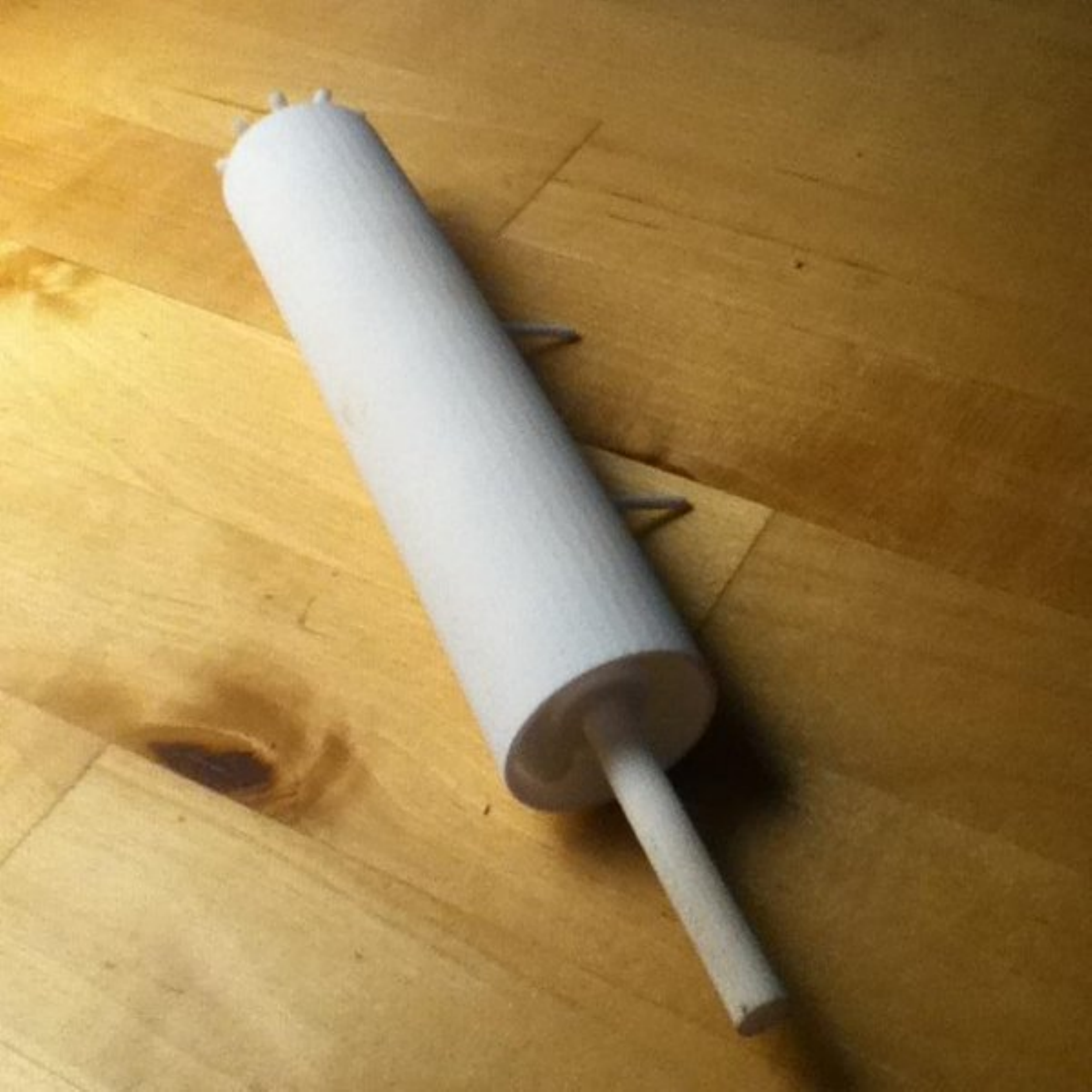}}
\caption{
The printed model can be assembled and used by
plugging it into a drilling machine.
}
\label{screw}
\end{figure}

\section{3D printing technology}

3D printing technology is a rapid prototyping process that allows the building of
objects using various technologies. The industry is relatively young. Only 30 years ago, 
the first patents of ballistic particle manufacturing were filed \cite{GRS,CLL}.
Jeremy Rifkin \cite{Rifkin,worldfinancialreview,economist2012} considers the process 
part of the ``third industrial revolution" which he describes as a time when 
``manufacturing becomes digital, personal, and affordable." \\

The first commercially successful 3D-printing technology appeared in 1994 with printed wax material.
More recently, with the technology of deposition of acrylate photopolymers,
the costs have dropped into a range of consumer technology. Printing services now 3D print in color
and high quality. 3D printing belongs to a larger class of construction methods called
rapid prototyping \cite{Cooper,CLL,KamraniNasr}, which started in the 1980s and include
a larger class of additive manufacturing processes. It is related
to other technologies like automated fabrication or computer numerical controlled 
machining \cite{GRS}. The development of 3D printing belongs to an ongoing chain
or evolution in industrialization and information. It is a trend created by 
different developments merging. Production, information technology, communication, and transportation are all 
linked closer than ever \cite{CLL}. \\

\begin{table}
\begin{tabular}{|l|l|l|l|}   \hline
{\bf Information revolutions}  & &   {\bf Industrial revolutions}     &                  \\  \hline
Gutenberg Press                & 1439    &   Steam Engine, Steal, Textile     &  1775    \\
Electric programmable computer & 1943    &   Automotive, Chemistry            &  1850    \\
Personal computer, Cell phone  & 1973    &   Laser printer, Rapid prototyping &  1969    \\
2G cellular technology         & 1991    &   3D printer                       &  1988    \\  \hline
\end{tabular} 

\caption{Information and industrial revolutions}
\end{table}

The attribute ``additive" in the description of 3D technology is used because the models  are
built by adding material. In sculpting or edging techniques by contrast, material 
is removed by chemical or mechanical methods.
3D printing is more flexible than molding \cite{Brain}. 
In the simplest case, the model is built 
from a series of layers, which ultimately form the cross sections of the object. 
Other 3D technologies allow to add material also from more than one direction. \\

In the last ten years, the technology has exploded and printers are now available 
for hobbyists. We have experimented with our own ``Up! StartPlus" printer
which is small enough to be carried into a calculus classroom for
demonstration. It allowed us for example to illustrate the concept of Riemann integral by 
illustrating the layer-by-layer slicing of a solid. \\

3D printing is now heavily used in the medical industry; the airplane industry; 
to prototype robots, art, and jewelry; nano structures; bicycles; ships; circuits;
make copies of artworks; or even to print houses or decorate chocolate cakes. 
Last and not least, and this is the thesis project of the second author, 
it is now also more and more used in education. 
While many courses deal primarily with the technique and design for 3D printing, 
a thesis of the second author focus on its use in mathematics 
education \cite{Slavkovsky}.

\section{Illustrations}

Archimedes is credited for the first attempts to find the area and circumference of a circle. 
He uses an exhaustion argument using the constant $\pi$ defined as the circumference divided
the diameter. To compute the area, he would compute the area of an inscribed regular polygon 
with contains the circle. Each polygon is a union of $n$ congruent triangles
whose union exhausts the disk more and more. To get an upper bound, he would look at polygons
containing the circle. The value of the circumference of the circle is so sandwiched between
the circumferences of the inner and outer polygon which converge. Modern calculus 
shoots all this down with explicit ``squeeze formulas" for the 
circumference $C(r)$ and area $A(r)$, which are
$C_i(r) = n r \sin(\pi/n) \leq C(r) \leq nr \tan(\pi/n) = C_o(r)$ and 
$rC_i(r) \leq A(r) \leq r C_o(r)$. The calculus student of today sees that Archimedes
was essentially battling the limit $\sin(x)/x$ for $x \to 0$. That the limit is $1$
is sometimes called the ``fundamental theorem of trigonometry" because one can deduce from this fact -
using addition formulas - the formulas the derivatives of the trig functions in general. We mention
the modern point of view because it shows that Archimedes worked at the core of 
calculus. He battled the notion of limit by comparing lower and upper bounds. 
The exhaustion argument also works for spheres, where it can be used that the 
volume $V(r)$ and surface area $A(r)$ of a sphere are related by $V=A r/3$. We illustrated
this with a 3D model. 

\begin{figure}
\scalebox{0.30}{\includegraphics{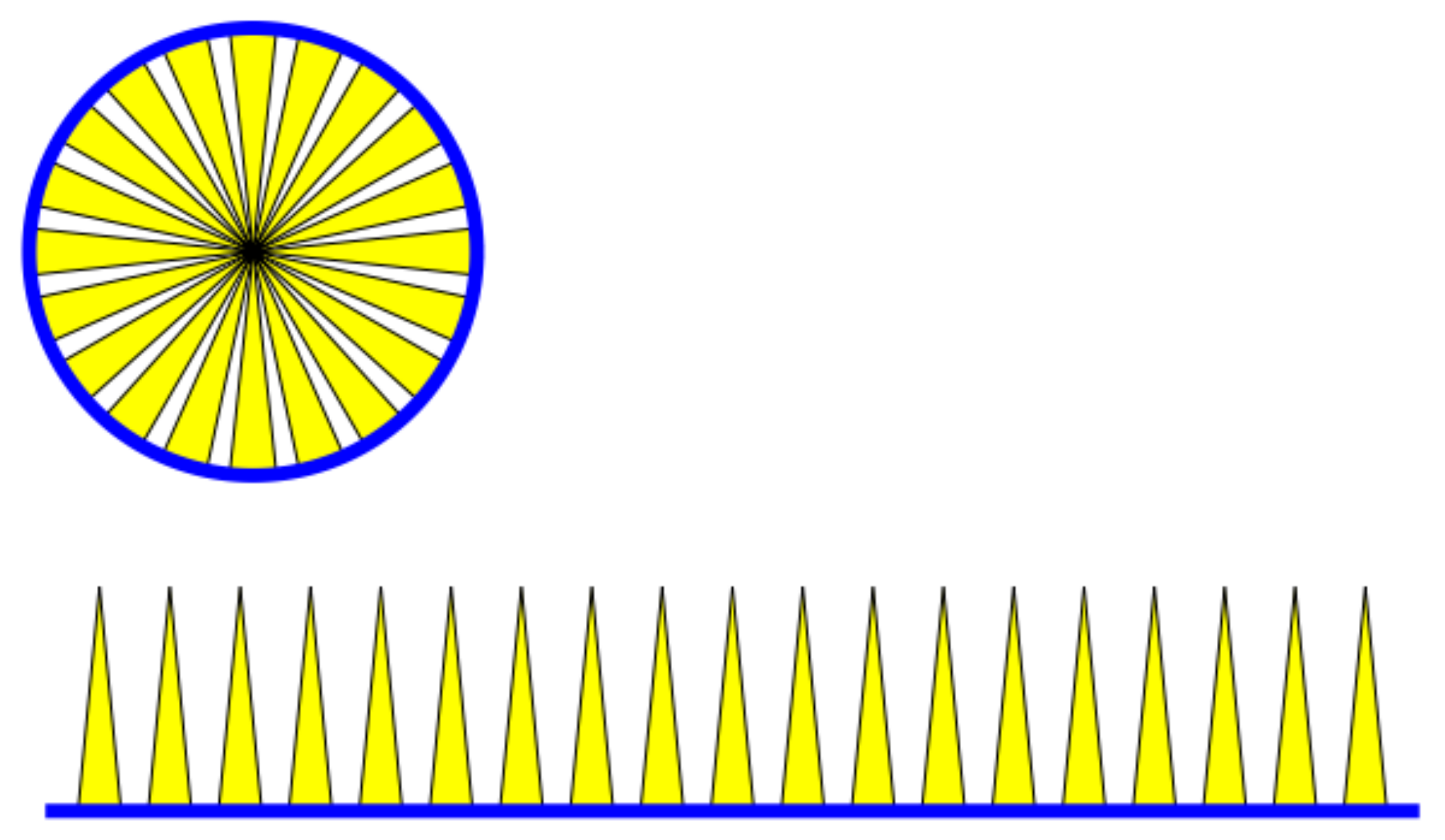}}
\scalebox{0.25}{\includegraphics{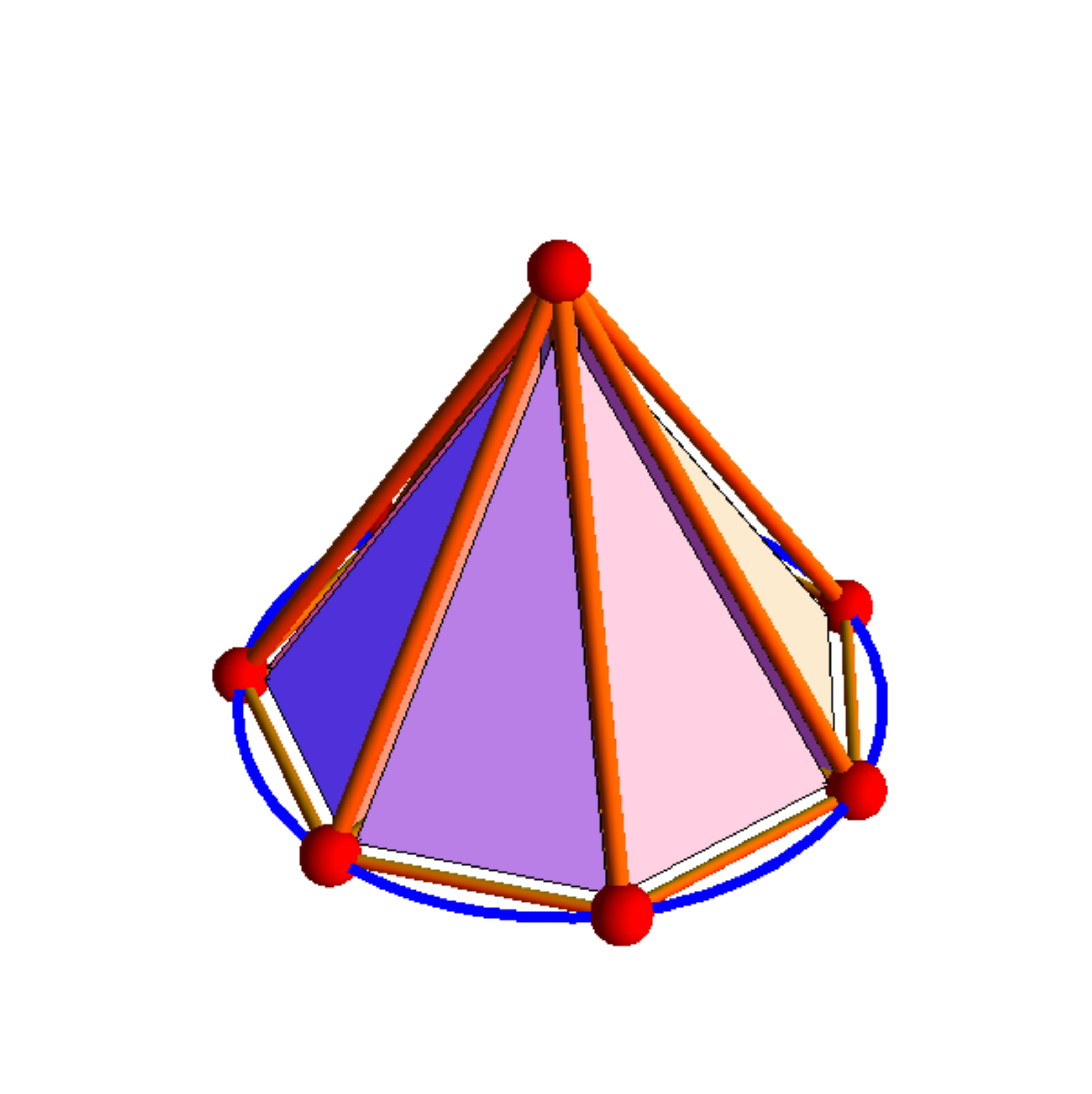}}
\caption{
The area of the circle as half the circumference times the radius is illustrated
by this figure of Archimedes. The argument can then be used to get the volume 
of a cone as $V=h A/3$ because each individual triangle of area $A$ 
of the triangularization of the circle produces a tetrahedron of volume $A h/3$. 
}
\label{archimedes}
\end{figure}

\begin{figure}
\scalebox{0.5}{\includegraphics{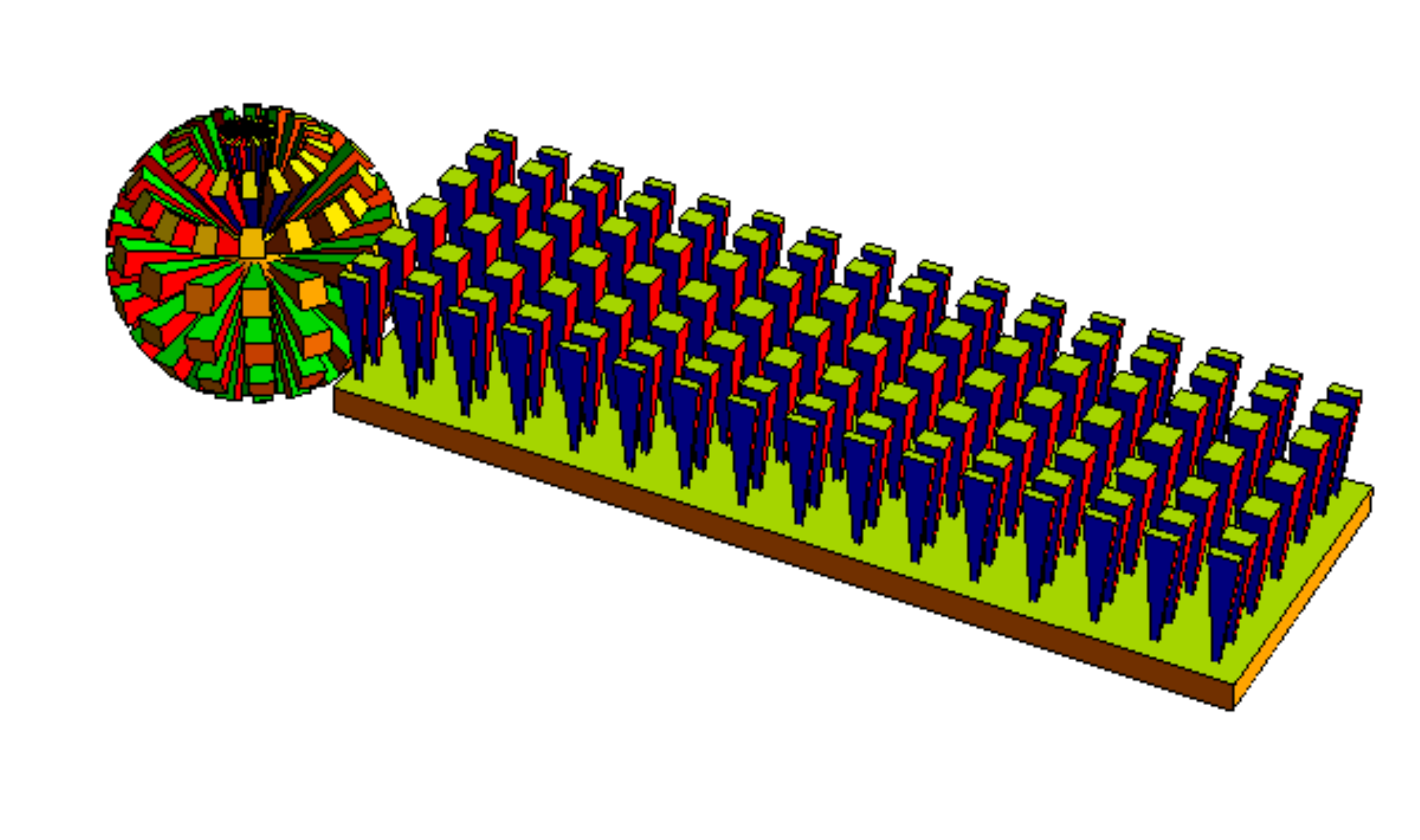}}
\caption{
We propose an Archimedes-style proof computing the volume of the sphere,
assuming that the surface area $A$ is known. The formula $V=A r/3$ can be 
seen by cutting up the sphere into many small tetrahedra of volume $dA r/3$.
When summing this over the sphere, we get $A r/3$.
}
\label{sphereproof}
\end{figure}

Archimedes' idea of comparing areas of cross sections of two
different bodies is today also called the {\bf Cavalieri principle} \cite{Kramer}.
In calculus classrooms, it is related to the ``washer method": if we slice a body $G$ 
perpendicularly along a line and if the area of a slice
is $A(z)$ and we use a coordinate system, where $z=0$ is the floor and 
$z=1$ the roof,  
$$   \int_0^1 A(z) \; dz $$ 
is the volume of the solid. Consequently, if two bodies can be sliced so that the 
area functions $A(z)$ agree, then their volumes agree. Archimedes understood the 
integral as a sum. Indeed, the Riemann sum with equal spacing 
$$  \frac{1}{n} \sum_{k=1}^n A(\frac{k}{n}) $$  
is sometimes called {\bf Archimedes sum} \cite{Agnew}. For a piecewise continuous
cross-section area $A(z)$, it is equivalent to the Riemann sum. 
The Cavalieri principle is a fundamental idea that perpetuates to other ideas of
mathematics.  \\

Of ``Archimedean style" \cite{Turnbull}, for example,
is the Pappus Centroid theorem, which determines the volume of a solid generated by rotating
a two-dimensional region in the $xz$ plane around the $z$ axes. Rotating two regions 
of the same area around produces the same volume. A mechanical proof is given in \cite{Levi}.
More general formulas can be obtained by taking a tubular neighborhood of a 
curve in space or a tubular neighborhood of a surface. The book \cite{Gray} 
has made it an important tool in differential geometry. \\

Greek mathematicians like Archimedes knew that the volume of a tetrahedron is one third 
of the product of the area of the base times the height \cite{Euclid}. One can see this by 
cutting a parallelepiped with twice the base area into six pieces. \\

This fact can be difficult to visualize on the blackboard. 
We used a 3D printer to produce pieces that add up to a cube.  \\

\begin{figure}
\scalebox{0.25}{\includegraphics{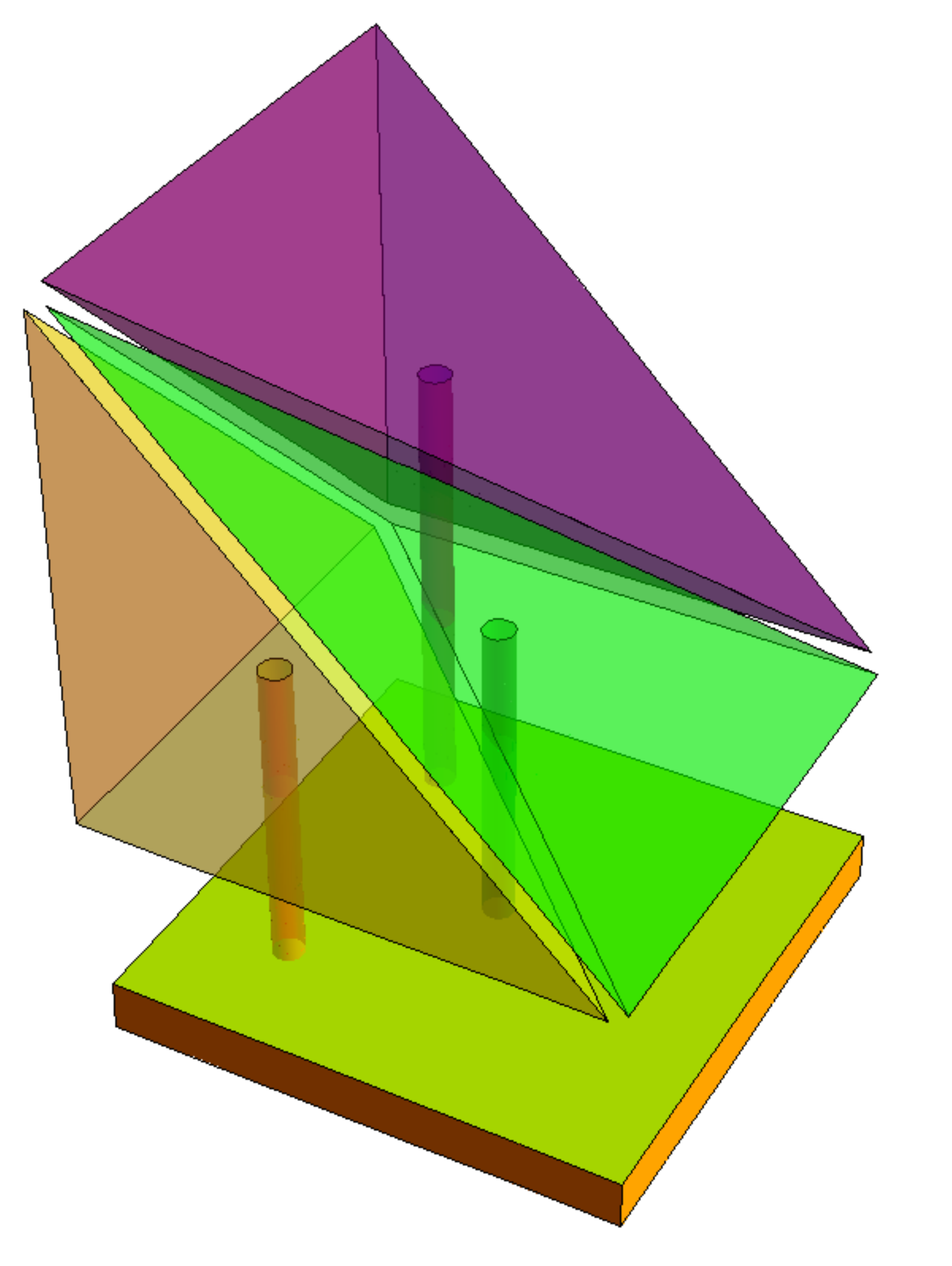}}
\scalebox{0.25}{\includegraphics{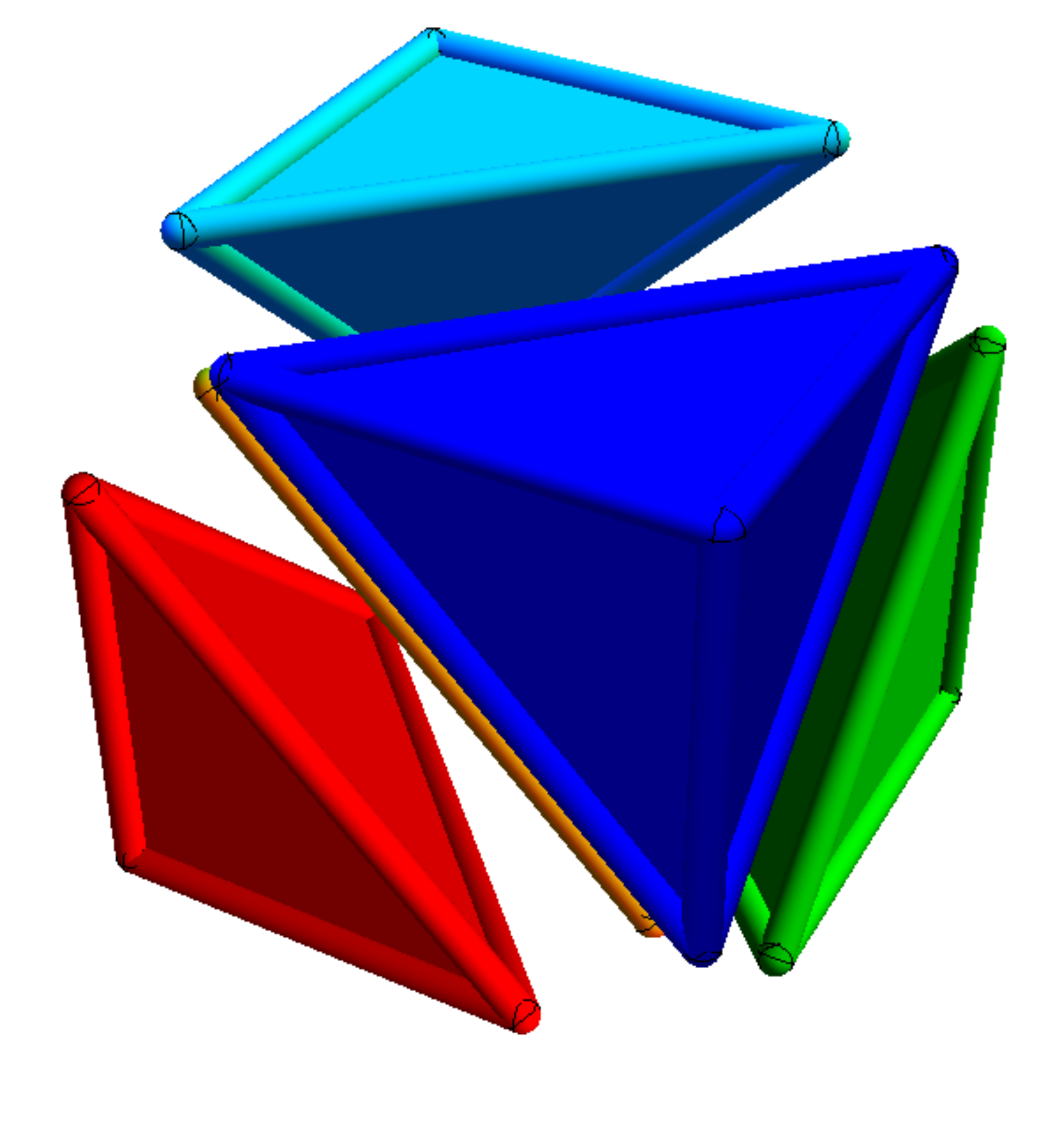}}
\caption{
Two attempts to visualize cutting a parallelepiped into 6 pieces of equal volume.
}
\label{parallelepieped}
\end{figure}

Given a two-dimensional region $G$ in the $xy$ plane and a point $P$
the convex hull defines a pyramid. The cross section at height $z$ has
area  $A(z) = A(0)(1 - (z/h)^2)$. Integrating this from $0$ to $h$ gives $A(0) h/3$.
To see this geometrically without calculus, one can triangulate the region
and connect each triangle $D$ to $P$. This produces tetrahedra which each have
volume $D h/3$. Summing over all triangles gives the approximate area of the base and 
summing all the tetrahedron volumes leads to the volume of the cone. \\

\begin{figure}
\scalebox{0.25}{\includegraphics{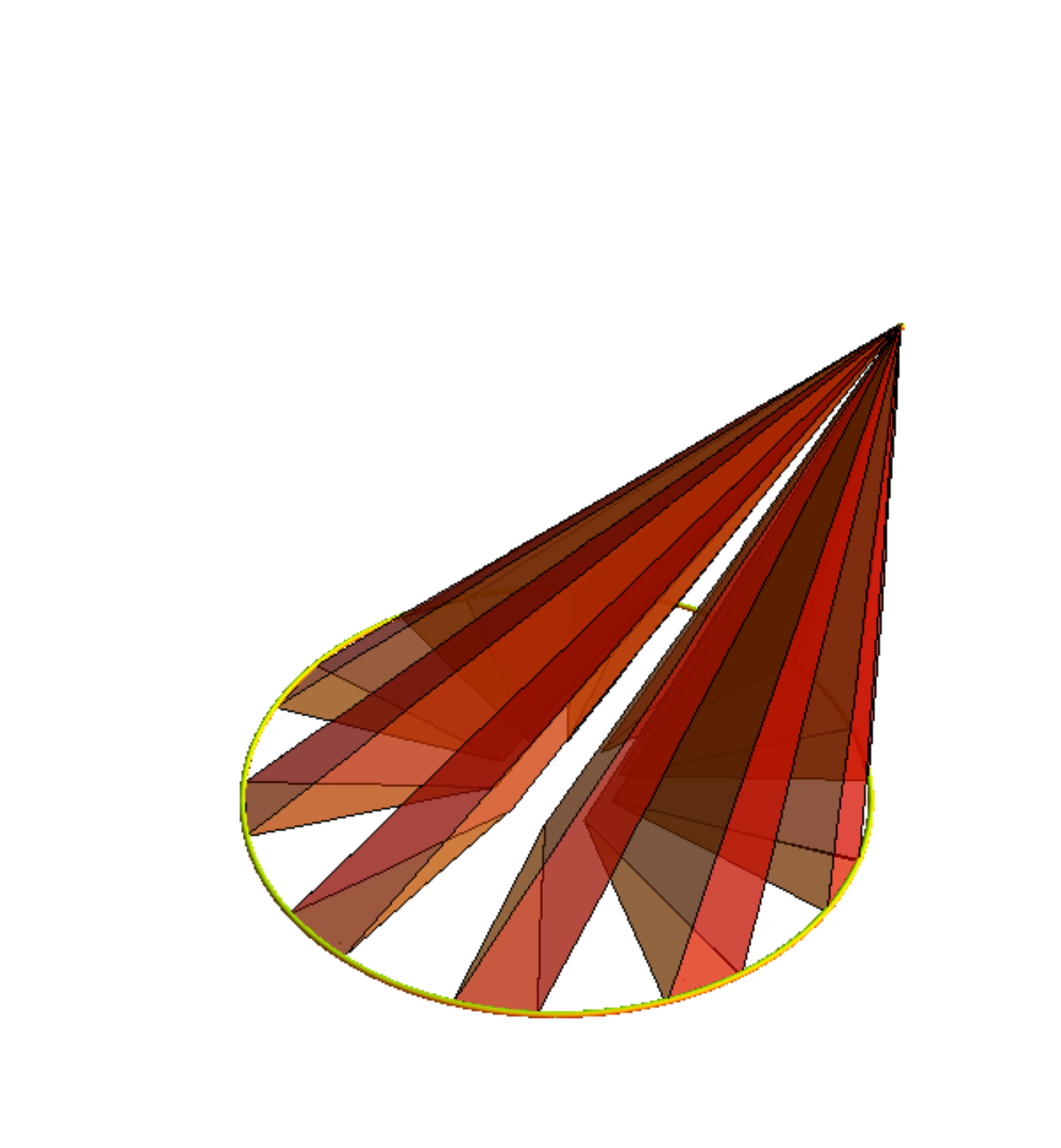}}
\scalebox{0.25}{\includegraphics{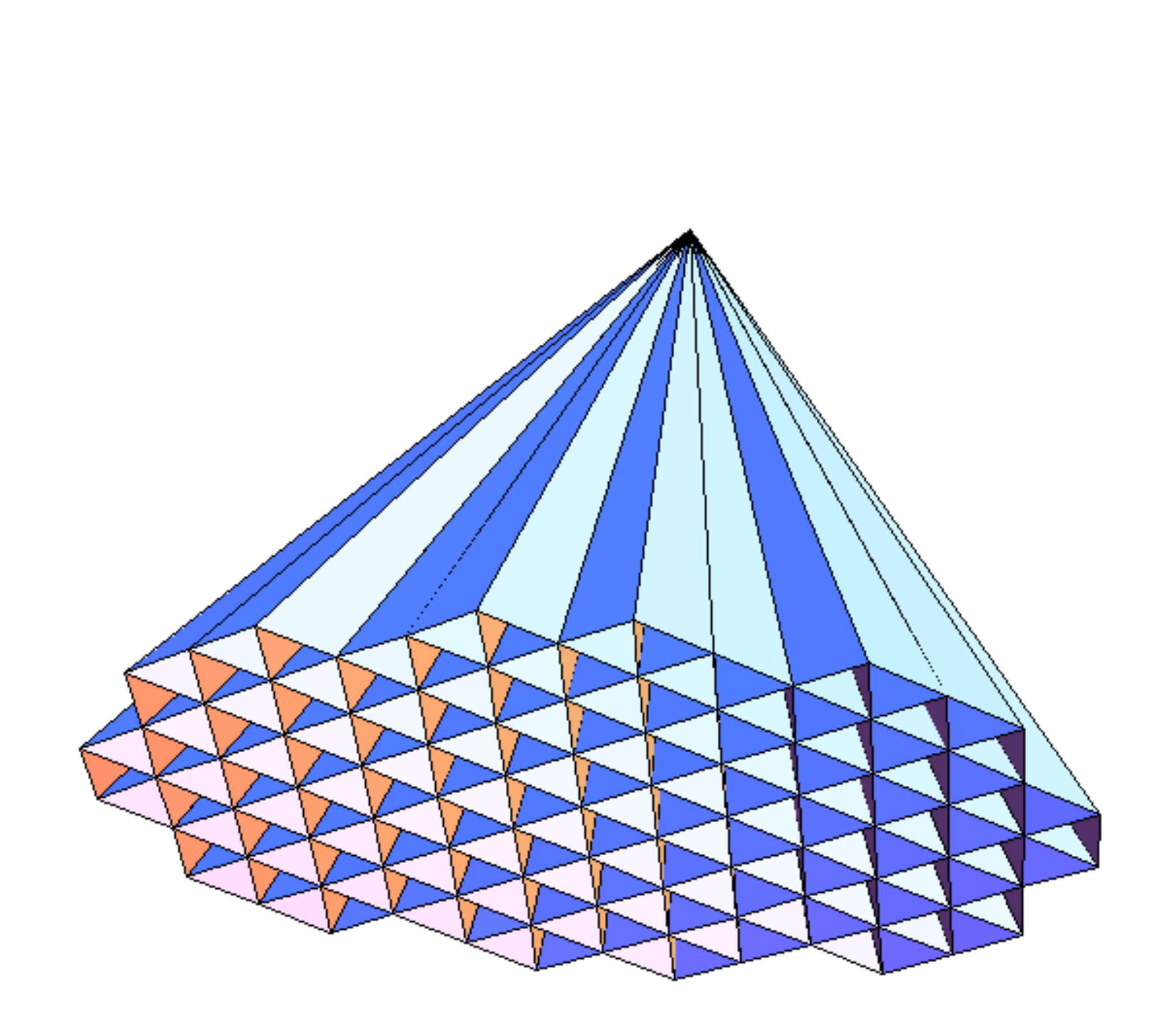}}
\caption{
A cone cut into small pieces reduces the volume computation of the cone to the 
volume computation of the pyramid. The later can be computed by cutting a cube into
6 pieces of equal volume. The volume formula $V=Ah/3$ for the pyramid with square
base was already known to the Greeks. The dissection method verifies it for general
cones. 
}
\label{archimedes}
\end{figure}

An impressive example of the comparison method is the computation
of the volume of the sphere. Archimedes used the same principle to compute the volume 
of Archimedian globes or domes (half globes) \cite{Apostol}. 
Archimedes compared the volume of a sphere with the differences of a cylinder
and a cone. For the sphere, the area of a cross section  is a disc of
radius $\sqrt{1-z^2}$ which leads to an area $\pi (1-z^2)$. This is also the 
area of the ring obtained by cutting the complement of a cone in a cylinder.

To illustrate this idea, we printed a half-sphere bowl that drains into 
the complement of a cone in a cylinder (see Figure~\ref{sphereproof}).
To demonstrate the proof, we fill in the top with water and let it drain 
into the base. We can now drink the proof. 
A 3D model that illustrates this can be seen in Figure~\ref{sphereproof}. In the case of the Archimedean 
dome, the comparison is done with the complement of a polygonal pyramid inside a circumscribing
prism verifies also here that the volume is $2/3$ of the prism. In the limit, when the 
polygon approaches a circle, one gets again that the volume of the sphere is $2/3$ times
the volume of the cylinder in which the sphere is inscribed. 

\begin{figure}
\scalebox{0.25}{\includegraphics{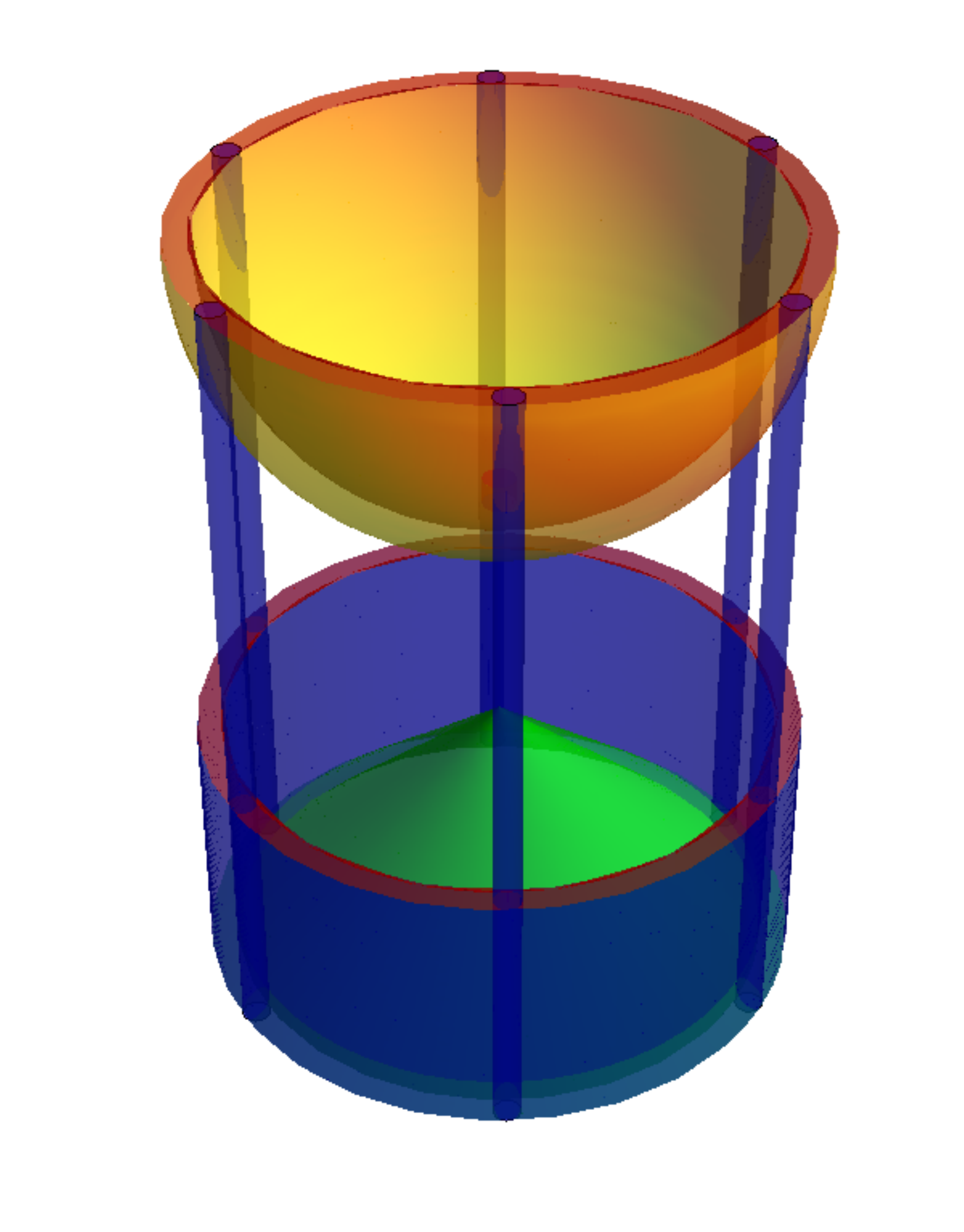}}
\scalebox{0.25}{\includegraphics{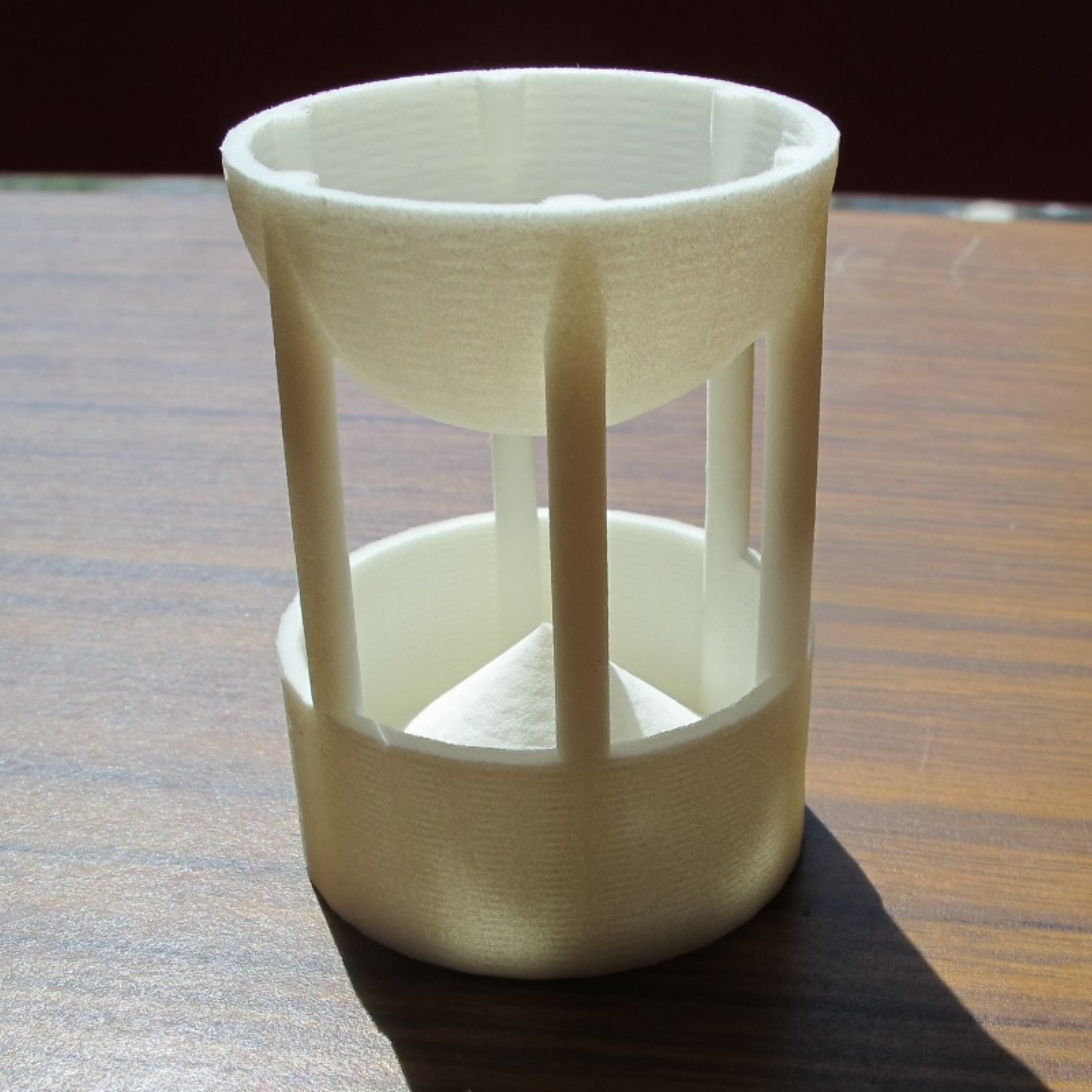}}
\caption{
The drinkable Archimedes Proof. The Mathematica model and the 
printout. As a demonstration, one can fill the spherical reservoir on top
with water. After it has dripped down, it fills the complement of a cone
in a cylinder. The volumes match. 
}
\label{archimedes}
\end{figure}

A classical problem asks for the volume of the intersection of two cylinders
intersecting perpendicularly leading to a ``Himmelsgew\"olbe". 
Also this object has been studied already by Archimedes \cite{Heath2}
and is an example of an Archimedean globe. 
The problem to compute its volume or the surface area appears in every calculus
book.  Another problem is the intersection of three cylinders, 
an integration problem that is also abundant in calculus books. 

\begin{figure}
\scalebox{0.20}{\includegraphics{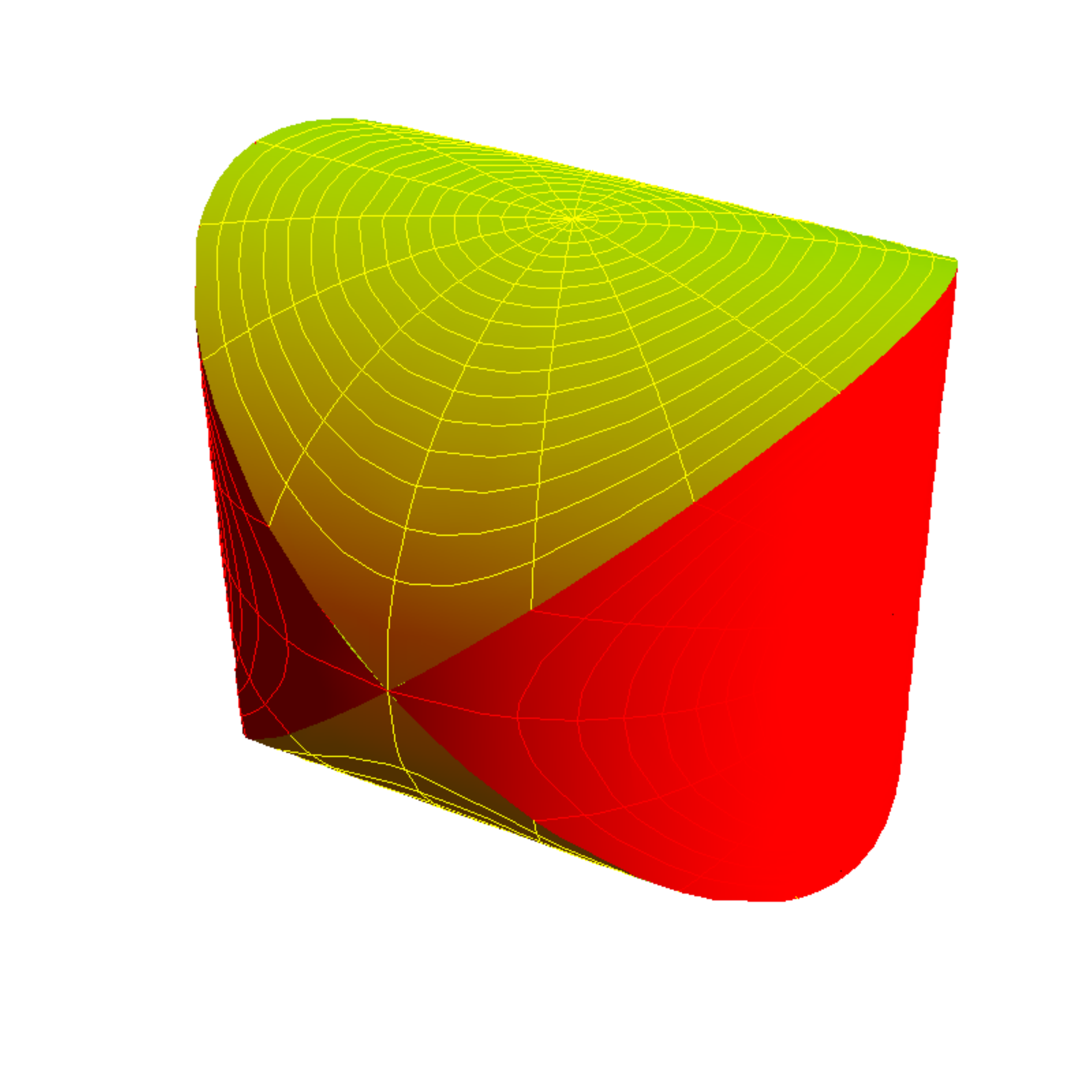}}
\scalebox{0.20}{\includegraphics{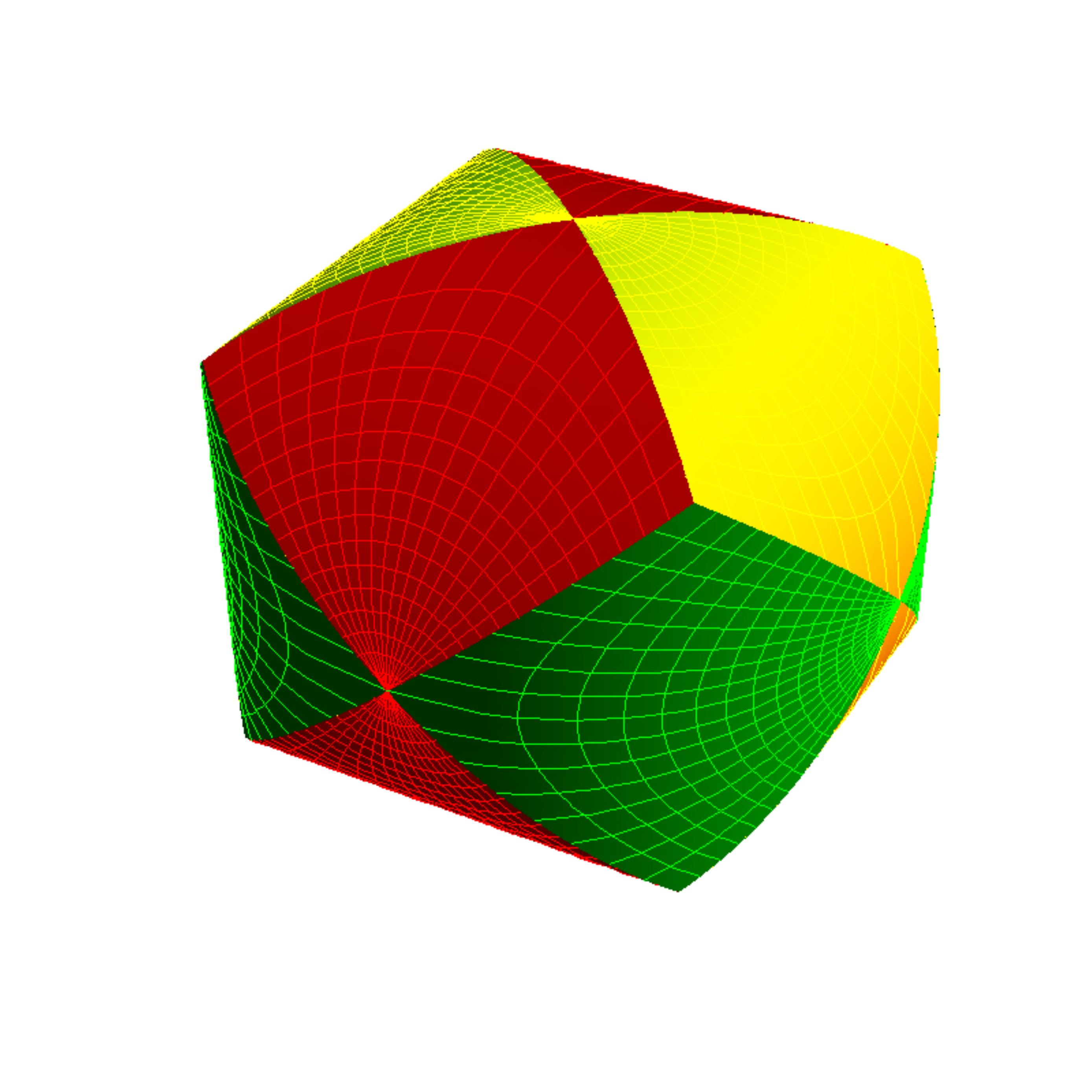}}
\caption{
Archimedes studied the solid obtained by intersecting 2 cylinders. 
This as well as the intersection of 3 cylinders is now a standard 
geometry problem in calculus textbooks. 
}
\label{screw}
\end{figure}

The first one is a solid which projects onto a square or two discs. 
The latter is a body which is not a sphere and projects on to three discs.  \\

Another solid that fits into this category is the 
{\bf Archimedian hoof}. It is the solid bound by a cylinder and two 
planes. The hoof plays an important role in 
the understanding of Archimedes thinking \cite{archimedescodex}.

\begin{figure}
\scalebox{0.20}{\includegraphics{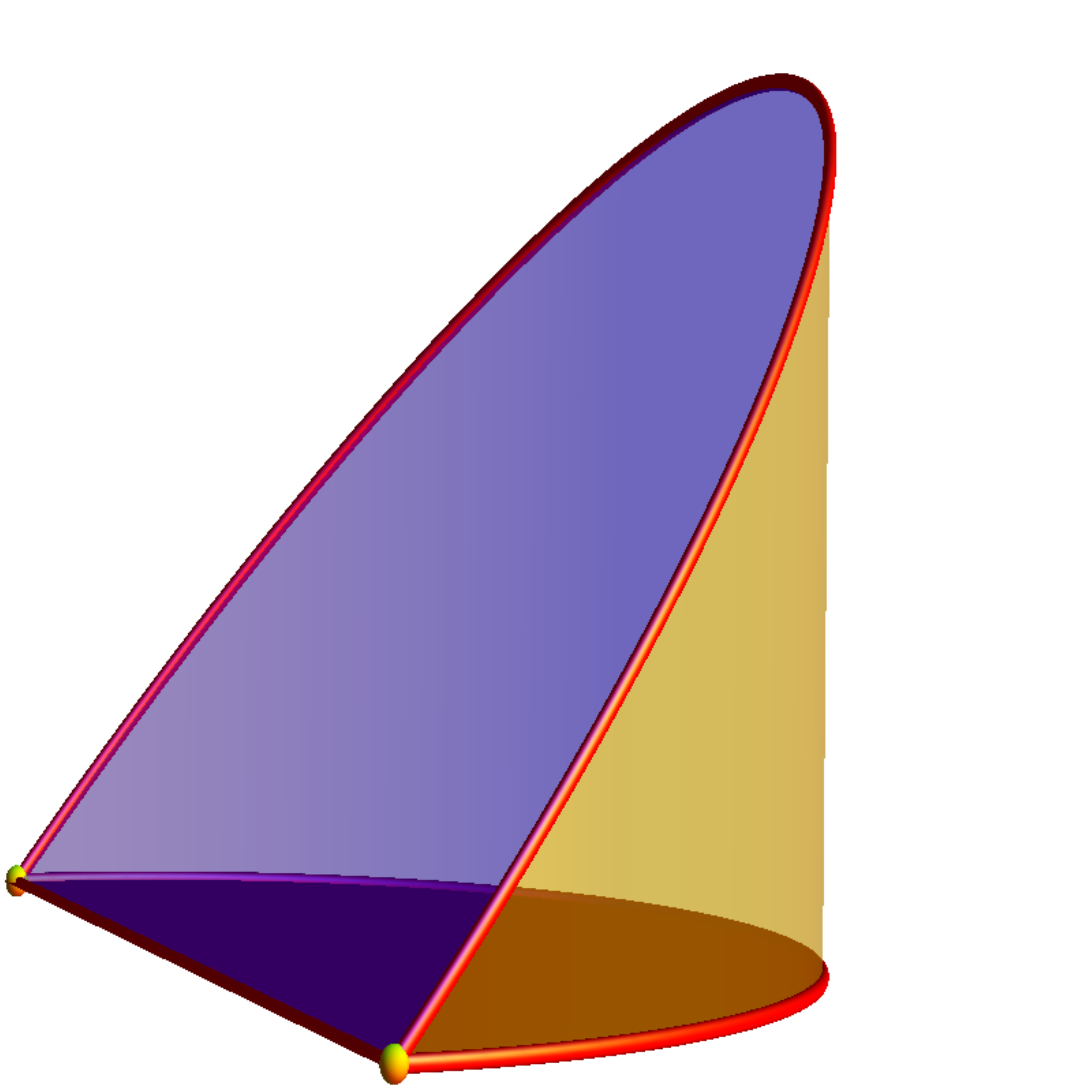}}
\scalebox{0.20}{\includegraphics{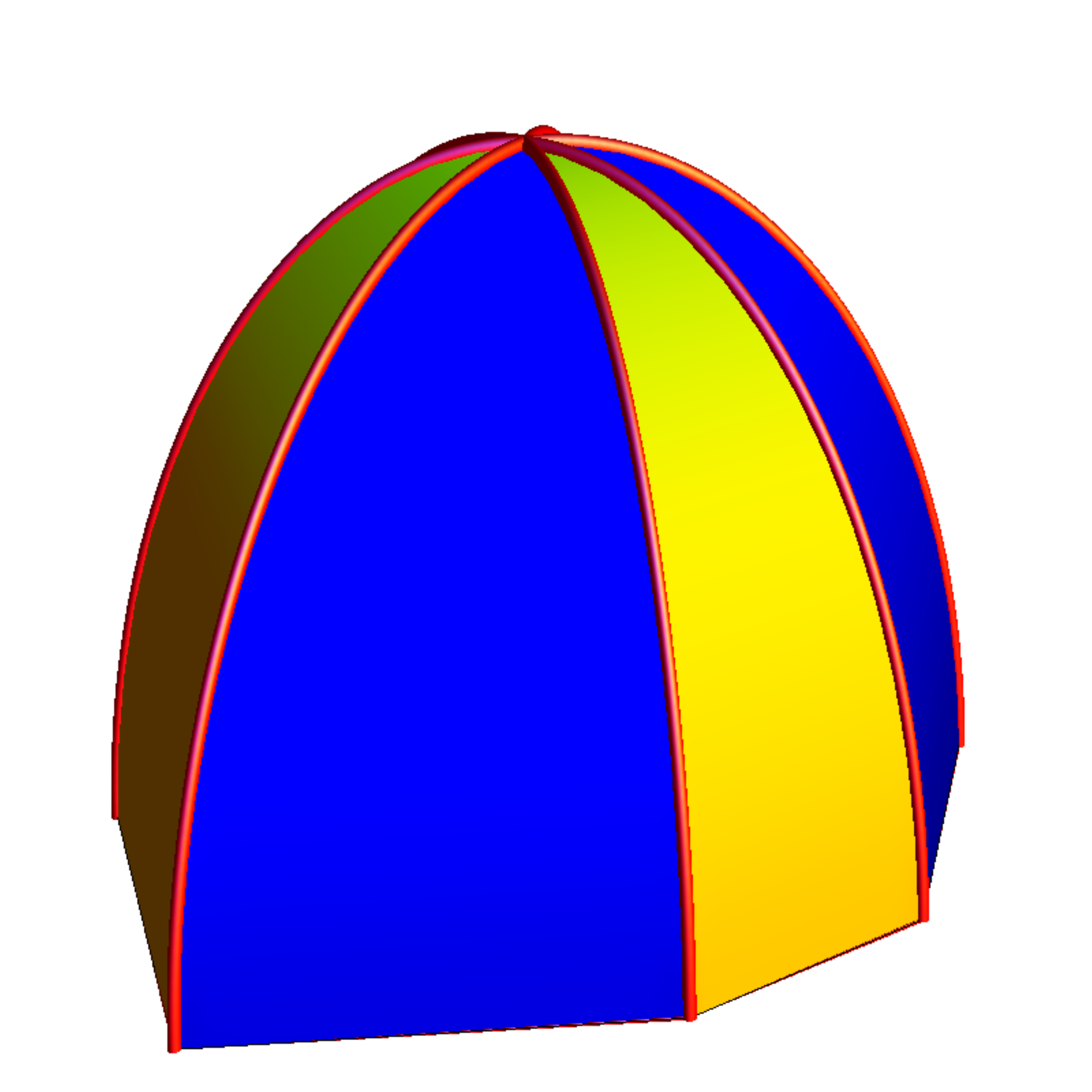}}
\caption{
The hoof of Archimedes and the Archimedean dome are solids for
which Archimedes could compute the volume with comparative integration
methods \cite{Apostol}. The hoof is also an object where Archimedes had to 
use a limiting sum, probably the first in the history of humankind
\cite{archimedescodex}.
}
\label{screw}
\end{figure}

The {\bf Archimedes hoof} \cite{Heath}
is the solid between the planes $z=x$ and $z=0$ inside the
cylinder $x^2+y^2 \leq 1$. It appears in any calculus textbook as an application in
polar integration
$$ \int_{-\pi/2}^{\pi/2} \int_0^1 r \cos(\theta) r \; dr d\theta  = \frac{4}{3} \; ,  $$
which Archimedes knew already to be $1/6$ of the volume $8$ of the cube
of side length $2$ containing the hoof. Archimedes figured out the volume differently.
He noticed that the slices $y=c$ lead to triangles which are all similar and have area $1-x^2$.
Integrating this up from $x=-1$ to $x=1$ is the area under the parabola $z=1-x^2$ and $z=0$,
a result which Archimedes already knew. The hoof actually might be the birth crib of 
infinitesimal calculus as we know it today. 

\section{Other bodies}

Martin Gardner mentions as part of nine problems in \cite{Gardner87} the question to 
find a volume of a body that fits snugly into a circle, square or triangle. 
Gardner hints that the volume is easy to determine without calculus. \\

Indeed, this problem is in the Archimedes style because the volume is half of the cylinder as
can be seen that when slicing it in the direction where one sees the triangle: at every 
cross section the triangle has half the area of the rectangle. Taking the rectangle
instead of the triangle produces the cylinder of area $\pi$. The cork therefore has volume 
$\pi/2$. 

\begin{figure}
\scalebox{0.25}{\includegraphics{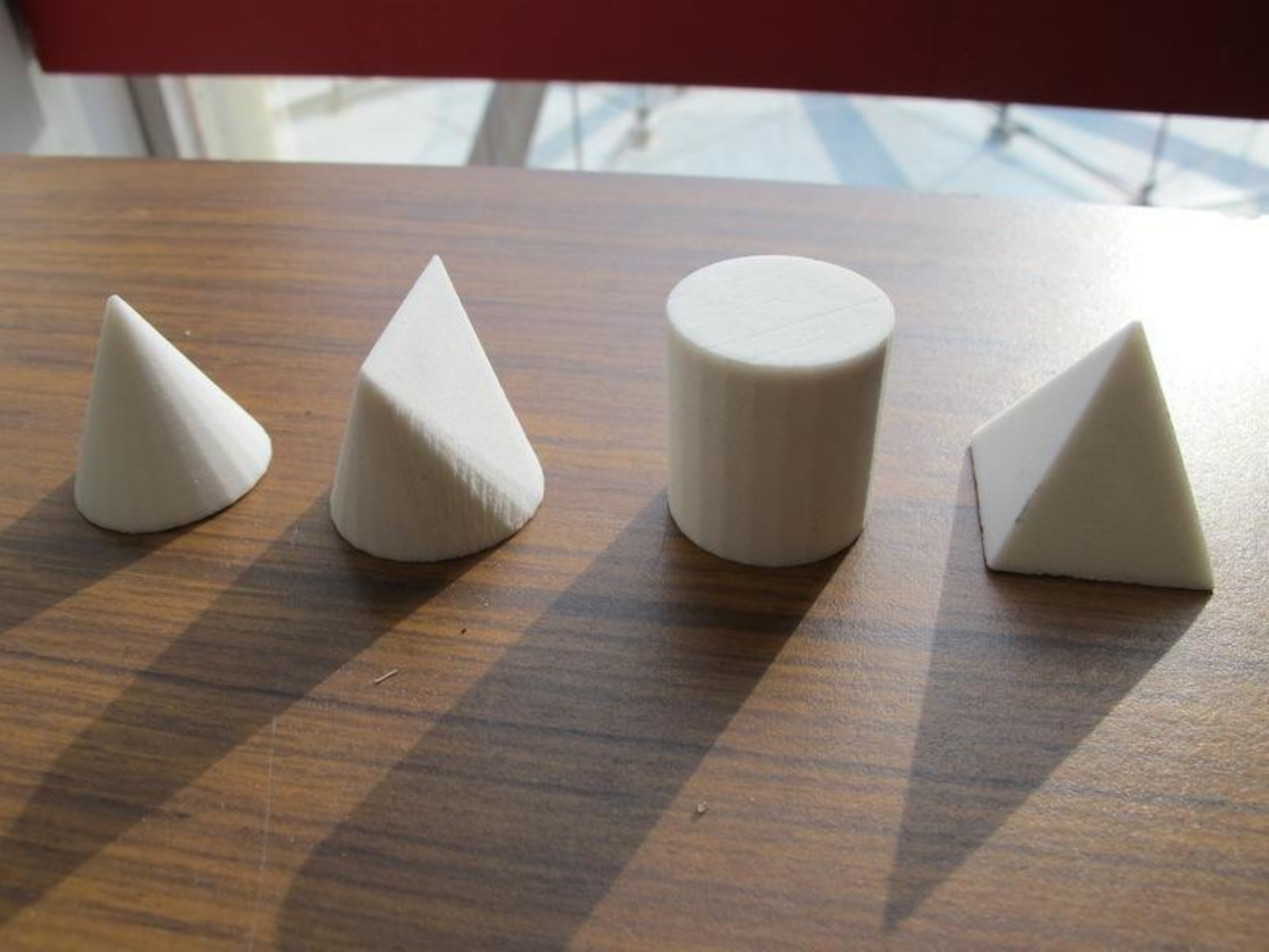}}
\scalebox{0.18}{\includegraphics{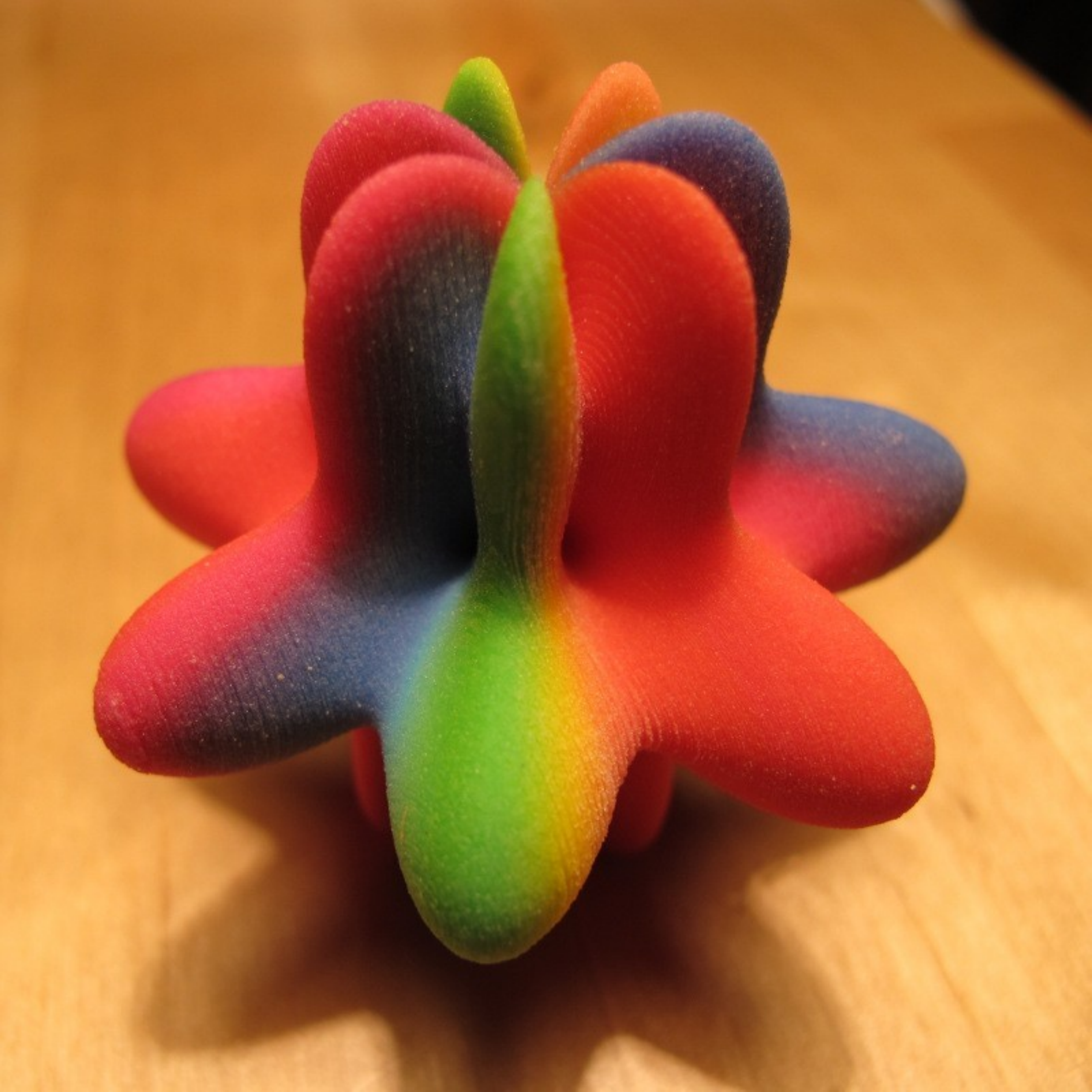}}
\caption{
Printing related to a problem of Martin Gardner motivated by 
\cite{Hart} who has used rapid prototyping to print this before. 
With printed version one can modify the problem in many ways like 
finding a solid where two cross sections are circles
and one cross section is an equilateral triangle, etc. To the right
we see a color printout of a surface built in Mathematica.
}
\label{cork}
\end{figure}

Finally, we illustrate another theme of circles and spheres which Archimedes would have 
liked but did not know. The Steiner chain and its three dimensional versions
illustrate a more sophisticated comparison method. 
Applying a M\"obius transformation in the flat case leads to chains where the spheres have 
different radius. If one extends the circles to spheres, one obtains Soddy's multiplets
\cite{Ogilvy}  These are shapes that are interesting to plot.  \\

\begin{figure}
\scalebox{0.20}{\includegraphics{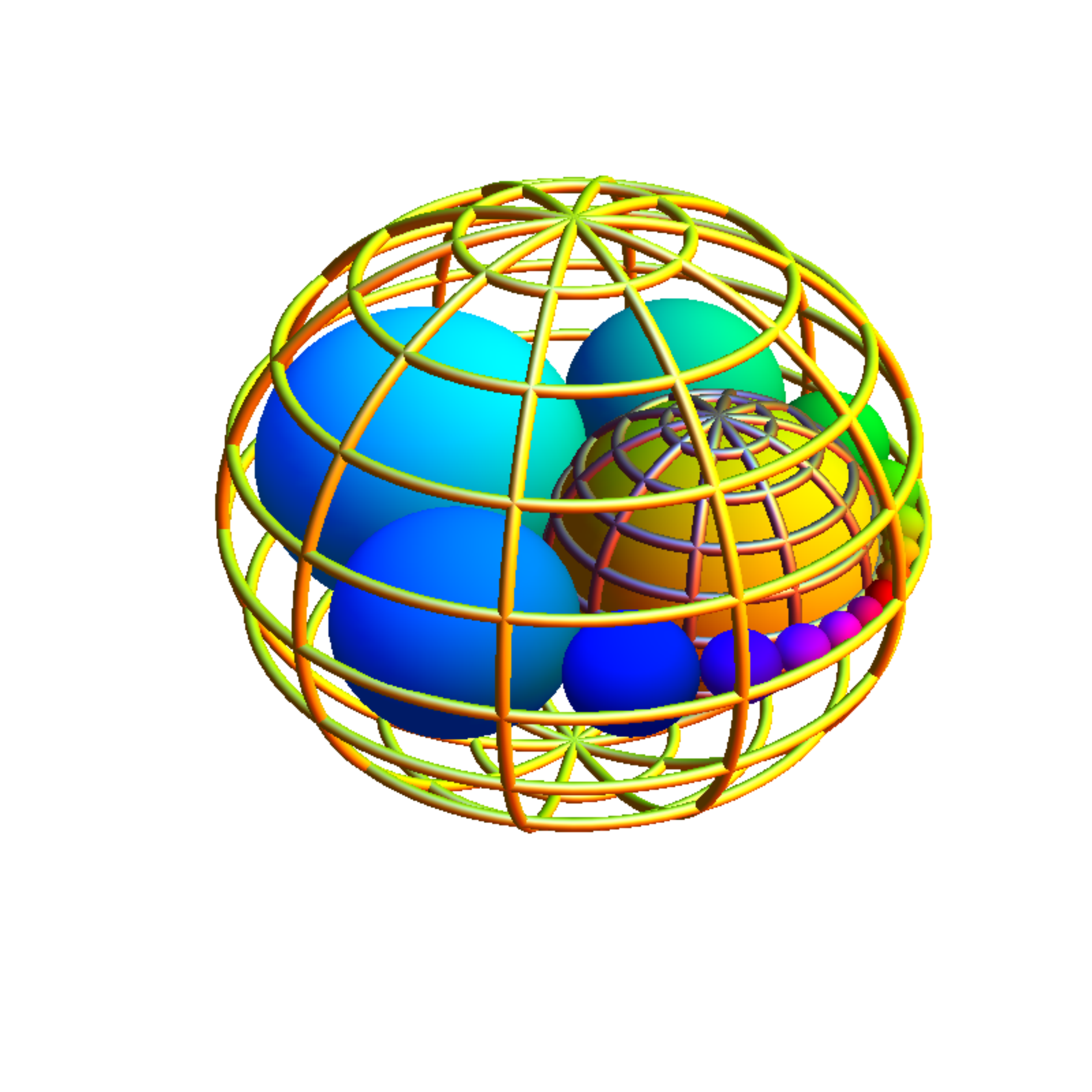}}
\caption{
Printing a three-dimensional Soddy multiplet. This solid illustrates
that 3D printers are capable of printing fairly complicated shapes. 
}
\label{screw}
\end{figure}

\begin{figure}
\scalebox{0.2}{\includegraphics{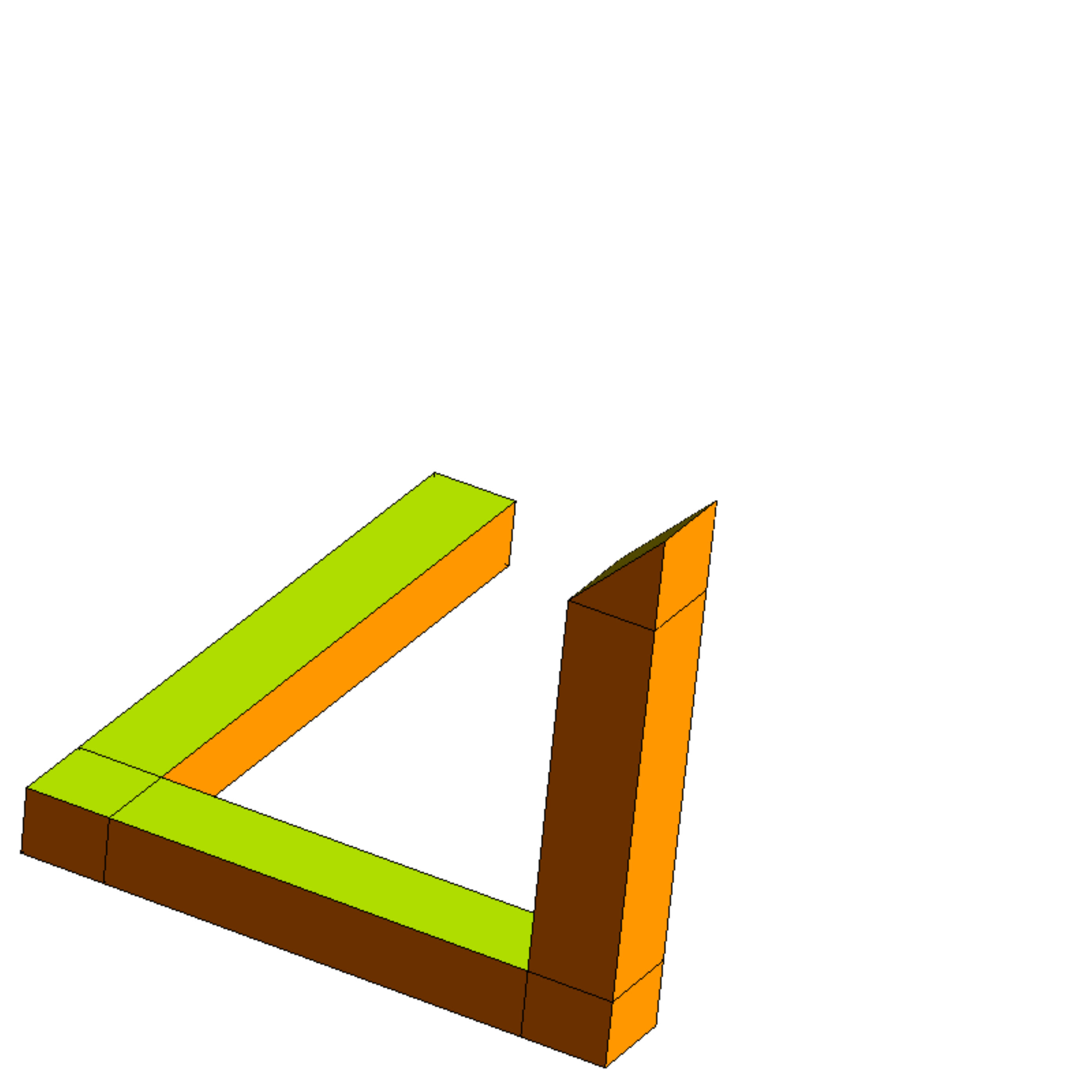}}
\scalebox{0.2}{\includegraphics{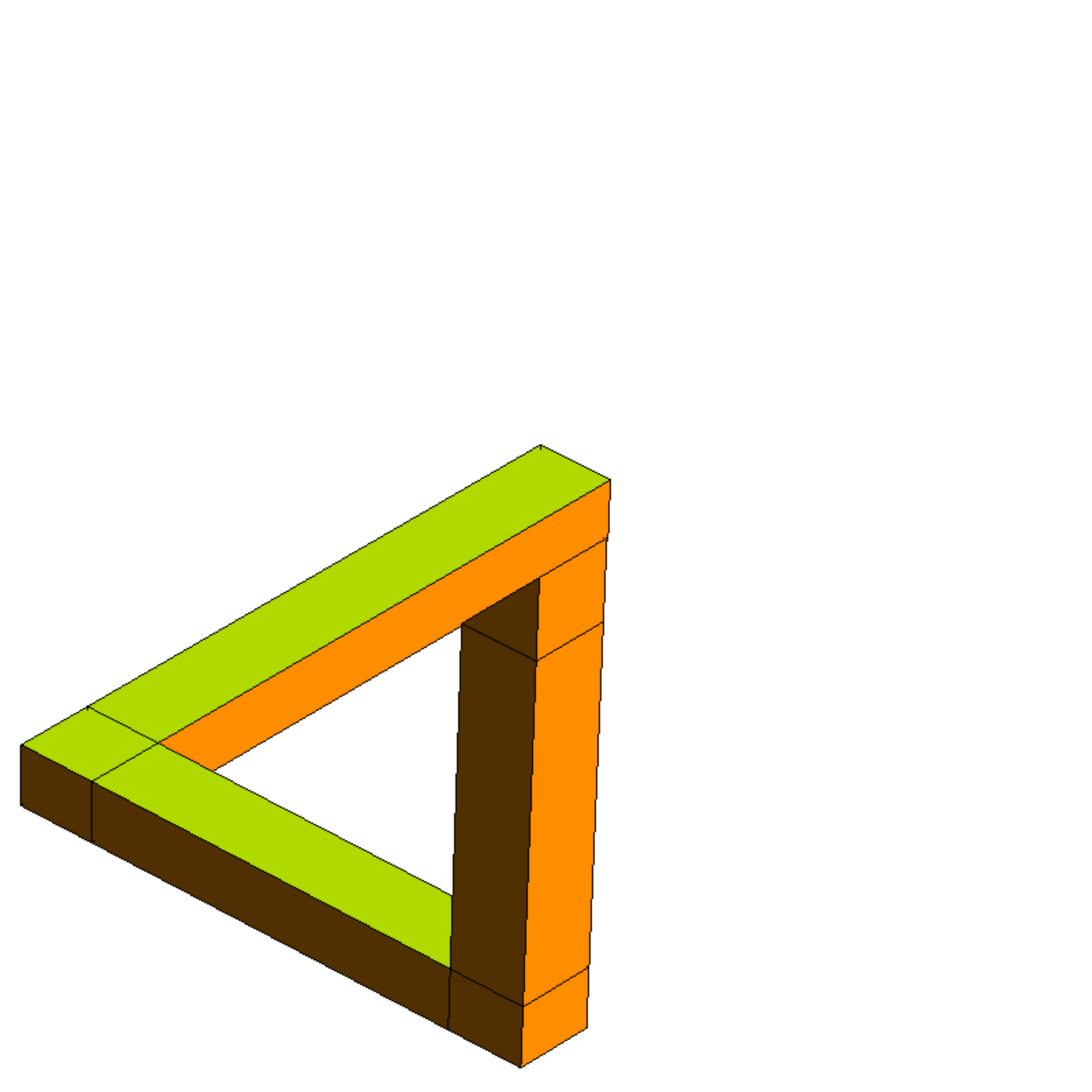}}
\caption{
Printing the Penrose triangle. The solid was created by Oscar Reutersvard and 
popularized by Roger Penrose. See \cite{Francis}. A Mathematica 
implementation has first appeared in \cite{Trott}. 
}
\label{archimedes}
\end{figure}

\section{3D printing and education}

We illustrated in this project also that this technology can lead to new perspectives
when looking at the history of mathematics and engineering. The new and still cutting-edge
technology allows everybody to build models for the classroom.  \\

Physical models are important for hands-on active learning. This has led to 
a repository of 3D printable models for education \cite{Lipson}. \\

The principles of integration appear in a natural way in 3D printing,
where a three dimensional object is built up automatically, like layer by layer.
Printing objects has become a fancy activity as printers have become more affordable.
One of us (E.S.)  purchased a 3D printer to experiment with
the technology to explore the technology for the classroom.  There are also online services
like ``Shapeways" or ``Sculpteo" that print an object and deliver it. 
Using 3D printing technology in education is hardly new
\cite{Lipson}. It has even been considered for \cite{sustainabledevelopement}
sustainable development in developing countries using sand in the desert. \\

3D printing technology has been used K-12 education in STEM projects \cite{rapman},
and elementary mathematics education \cite{contemporary}. There is optimism that it
will have a large impact in education \cite{Cliff}.
\cite{BullGroves} write:
``As fabrication tools become increasingly accessible, students will be 
able to learn about engineering design and experience the thrill of 
seeing their ideas realized in physical form". We believe the same is true in mathematics.  \\

In \cite{Wohler}, Timothy Jump, who teaches at Benilde St. Margaret's Highschool is quoted as, 
``3D printing stimulates a student's mechanical-spatial awareness in ways that textbooks cannot." 
The article concludes with, 
``The use of 3D printing technology in education is growing quickly because educators 
at all levels and from many areas of expertise can see how the technology helps students 
learn in very real, visual, and tangible ways." \\

The technology is fantastic to illustrate concepts in calculus. One of us (E.S.) 
has taught part of a lecture at Harvard using this technology. 
We have brought our printer into the classroom to illustrate the concept with 
single variable calculus students. \\

It has also been used in a multivariable calculus courses at Harvard as a final project.
Funding from the Elson Family Arts initiative is currently being used to print the 
student projects. 

\section{Tips for printing in 3D}

We conclude this article with some tips and remarks for anybody who wants to explore the subject more. 
First of all, the technology needs patience. Don't expect that things work immediately 
especially when printing on your own printer. You  need a lot of time and ability to overcome 
obstacles. Printer parts can fail, calibration can be subtle, software glitches occur frequently.
Commercially available printing services can refuse to print certain parts. 
It is still an area for enthusiasts. In addition, printers costing dozens of thousand 
dollars can break. We have seen this in local art shops or when submitting to 
professional printing services. It can happen for example that objects pass initial
tests but cannot be printed in the end. The sphere volume proof, for
example, was rejected several times, even after multiple optimization attempts. 
Another source for difficulties can be the size of the object, as it can affect the printing cost. \\

During the time of this project, the prices have fallen and the printers
have become easier to work with. This is good news for anybody who wants to experiment
with mathematical shapes beyond pen, paper, and computers.  \\

The technology of 3D printing is exciting because it develops fast. 
The prices have dropped to as little as 500 dollars.  
Some of the insights we gained when working with printing could be useful also in the future, when 
printers will be a standard in many households. 
Here are few tips when dealing with 3D printing in education: 

\begin{itemize}
\item To print more complicated objects or objects which have movable parts, 
add thin connections to be removed before use. An example is the 
Archimedes cochlias which contains a spiral moving inside a cylinder. 
\item  Most printers use millimeters for units, 
so it is best to design using a millimeter template. The same
units work best for importing into a printing program or printing service. 
\item It is important to check the tolerances of the 
printer and the materials before designing a shape. 
\item 
Walls cannot be too thin. Materials and printers require at least 1mm thickness, and sometimes up to 3mm, 
and this is another reason to design in millimeters and not inches. 
\item Connecting cylinders that provide stability cannot be too narrow or they break.
\item Hollow shapes are less expensive to print than solids.
\item Solid parts need to be closed to print correctly.
\item The nozzle and platform of the printer need to be heated to the correct temperature.
\item Ensure that the printout is in the correct STL format, either ASCII or binary. 
      Software like ``Meshlab'' or ``Netfabb" allows to work on the STL format
      (unfortunately, they both do not import WRL format). 
      Mathematica only exports in binary, whereas the UP! printer requires ASCII.
\item Not all modeling programs are equal. While Google SketchUp is easy to use, it is difficult 
      to make accurate models. Computer algebra systems are better suited in this respect. 
\item The design program should export STL or 3DS format. 
      ``SketchUp" requires a plugin to export to STL. 
      Mathematica can export both to STL and WRL, the later even in color. For some reason, only 
      Export[``file.wrl",Graphics] works but Export[``file.wrl",Graphics,``WRL"] does not. 
\item Printing in color works in principle, but not always. Since the export and import filters
      are still young, most color outputs do not print in color. 
\item The orientation can matter. To print a cylinder for example, 
      it can better to print along the axes. 
\item Printing yourself can take hours. Printing with printing services
      can need patience too. Typically we obtained the objects within a week.
\item There are a lot of materials already available. Plastic, 
      ceramics, sandstone, metal and even chocolate is possible.
      Different materials have different tolerances.
\item The price depends very much on the size of the model. 
\item Sometimes, the object is accepted by a printing service 
      and only discovers later that it is not possible. 
\item When printing in color from Mathematica, export in WRL 
      (Virtual {\bf W}orld {\bf R}eality Modeling {\bf L}anguage) format.
      (VRML currently does not print). 
\end{itemize}

\vspace{12pt}

\begin{thebibliography}{10}

\bibitem{Agnew}
R.~P. Agnew.
\newblock {\em Calculus, Analytic Geometry and Calculus, with Vectors}.
\newblock McGraw-Hill Book Company, inc, 1962.

\bibitem{Apostol}
T.M. Apostol and M.A. Mnatsakanian.
\newblock A fresh look at the method of {A}rchimedes.
\newblock {\em American Mathematical Monthly}, 111:496--508, 2004.

\bibitem{Archimedes}
Archimedes.
\newblock On the sphere and cylinder.
\newblock In {\em The works of Archimedes}. C.J. Clay and Sons, Cambridge
  University Press Warehouse, 1544 (reprint 1897).
\newblock Widener Math 275.1.105, 275.1.41.3.

\bibitem{Becker}
O.~Becker.
\newblock {\em Das Mathematische Denken der Antike}.
\newblock Goettingen, VandenHoeck und Ruprecht, 1957.

\bibitem{Brain}
M.~Brain.
\newblock How {Stereolithography} {3-D} layering works.
\newblock {\\}http://computer.howstuffworks.com/stereolith.htm/printable, 2012.

\bibitem{BullGroves}
G.~Bull and J.~Groves.
\newblock The democratization of production.
\newblock {\em Learning and Leading with Technology}, 37:36--37, 2009.

\bibitem{CLL}
C.S.~Lim C.K.~Chua, K.F.~Leong.
\newblock {\em Rapid Prototyping}.
\newblock World Scientific, second edition, 2003.

\bibitem{Cliff}
D.~Cliff, C.~O'Malley, and J.~Taylor.
\newblock Future issues in socio-technical change for uk education.
\newblock {\em Beyond Current Horizons}, pages 1--25, 2008.
\newblock Briefing paper.

\bibitem{Cooper}
K.G. Cooper.
\newblock {\em Rapid Prototyping Technology, Selection and Application}.
\newblock Marcel Dekker, Inc, 2001.

\bibitem{Czwalina}
A.~Czwalina.
\newblock {\em Archimedes}.
\newblock B.G. Teubner, Leipzig und Berlin, 1925.

\bibitem{dijksterhuis}
E.J. Dijksterhuis.
\newblock {\em Archimedes}.
\newblock Princeton University Press, 1987.

\bibitem{Dunham}
W.~Dunham.
\newblock {\em Journey through Genius}.
\newblock Wiley Science Editions, 1990.

\bibitem{economist2012}
Economist.
\newblock The third industrial revolution.
\newblock {\em Economist}, Apr 21, 2012, 2012.

\bibitem{Euclid}
T.Heath Euclid.
\newblock {\em The thirteen Books of the Elements, Volume 3}.
\newblock Dover, 1956.

\bibitem{Eves}
H.~Eves.
\newblock {\em Great moments in mathematics (I and II)}.
\newblock The Dolciani Mathematical Expositions. Mathematical Association of
  America, Washington, D.C., 1981.

\bibitem{Francis}
G.~Francis.
\newblock {\em A topological picture book}.
\newblock Springer Verlag, 2007.

\bibitem{Gardner87}
M.~Gardner.
\newblock {\em Book of Mathematical Puzzles and Diversions}.
\newblock University of {C}hicago Press, 1987.

\bibitem{Geymonat}
M.~Geymonat.
\newblock {\em The Great Archimedes}.
\newblock Baylor University Press, 2010.

\bibitem{Gray}
A.~Gray.
\newblock {\em Tubes}.
\newblock Addison-Wesley Publishing Company Advanced Book Program, Redwood
  City, CA, 1990.

\bibitem{Hart}
G.~Hart.
\newblock Geometric sculptures by {George Hart}.
\newblock http://www.georgehart.com.

\bibitem{Heath3}
T.L. Heath.
\newblock {\em A history of Greek Mathematics, Volume II, From Aristarchus to
  Diophantus}.
\newblock Dover, New York, 1981.

\bibitem{Heath}
T.L. Heath.
\newblock {\em The works of Archimedes}.
\newblock Dover, Mineola, New York, 2002.

\bibitem{Heath4}
T.L. Heath.
\newblock {\em A Manual of Greek Mathematics}.
\newblock Dover, 2003 (republished).

\bibitem{Heath2}
T.L. Heath.
\newblock {\em The Method of Archimedes, Recently discovered by Heiberg: A
  supplement to the Works of Archimedes}.
\newblock Cosimo Classics, 2007.

\bibitem{Heiberg}
J.L. Heiberg.
\newblock {\em Geometrical Solutions Derived from Mechanics, A Treatise of
  {A}rchimedes}.
\newblock The open court publishing Co, 1909 (translated 1942).

\bibitem{Hirshfeld}
A.~Hirshfeld.
\newblock {\em Eureka Man, the Life and Legacy of Archimedes}.
\newblock Walker and Company, New York, 2009.

\bibitem{GRS}
D.~Rosen I.~Gibson and B.~Stucker.
\newblock {\em Additive Manufacturing Technologies}.
\newblock Springer, 2010.

\bibitem{BBG}
R.~Girgensohn J.~Borwein, D.Bailey.
\newblock {\em Experimentation in Mathematics}.
\newblock A.K. Peters, 2004.
\newblock Computational Paths to Discovery.

\bibitem{Jaeger}
Mary Jaeger.
\newblock {\em Archimedes and Roman Imagination}.
\newblock University of Michigan Press, Ann Arbor, 2008.

\bibitem{KamraniNasr}
A.~Kamrani and E.A. Nasr.
\newblock {\em Rapid Prototyping, Theory and Practice}.
\newblock Springer Verlag, 2006.

\bibitem{Katz2011}
V.J. Katz.
\newblock {\em A history of Mathematics}.
\newblock Addison-Wesley, 1998.

\bibitem{Kramer}
E.E. Kramer.
\newblock {\em The Nature and Growth of Modern Mathematics, Volume 1}.
\newblock Fawcett Premier Book, 1970.

\bibitem{Kropp}
G.~Kropp.
\newblock {\em Geschichte der Mathematik}.
\newblock Quelle und Meyer, 1969.

\bibitem{rapman}
G.~Lacey.
\newblock 3d printing brings designs to life.
\newblock {\em techdirections.com}, 70 (2):17--19, 2010.

\bibitem{Levi}
M.~Levi.
\newblock {\em The Mathematical Mechanic}.
\newblock Princeton University Press, 2009.

\bibitem{Lipson}
H.~Lipson.
\newblock Printable 3d models for customized hands-on education.
\newblock Paper presented at Mass Customization and Personalization (MCPC)
  2007, Cambridge, Massachusetts, United States of America, 2007.

\bibitem{Lurje}
S.J. Lurje.
\newblock {\em Archimedes}.
\newblock Dr. Walter Hollitscher, Phoenix Buecherei, 1948.

\bibitem{lacuna}
R.~Masi{\`a}-Fornos.
\newblock A ``lacuna'' in {P}roposition 9 of {A}rchimedes' {\it {o}n the sphere
  and the cylinder}, {B}ook {I}.
\newblock {\em Historia Math.}, 37(4):568--578, 2010.

\bibitem{Netz}
R.~Netz.
\newblock {\em The works of Archimedes}.
\newblock Cambridge University Text, 2004.

\bibitem{archimedescodex}
R.~Netz and W.~Noel.
\newblock {\em The {A}rchimedes Codex}.
\newblock Da Capo Press, 2007.

\bibitem{Ogilvy}
C.S. Ogilvy.
\newblock {\em Excursions in Geometry}.
\newblock Oxford University Press, New York, 1969.

\bibitem{archimedes2}
Domenico~Fetti (Painter).
\newblock Archimedes portrait.
\newblock Dresden,
  {http://archimedes2.mpiwg-berlin.mpg.de/archimedes$\_$templates/popup.htm},
  1620.

\bibitem{sustainabledevelopement}
J.M. Pearce, C.M. Blair, K.J. Kaciak, R.Andrews, A.~Nosrat, and
  I.~Zelenika-Zovko.
\newblock 3-d printing of open source appropriate technologies for
  self-directed sustainable development.
\newblock {\em Journal of Sustainable Development}, 3 (4):17--28, 2010.

\bibitem{Pippenger}
N.~Pippenger.
\newblock Two extensions of results of {A}rchimedes.
\newblock {\em American Mathematical Monthly}, 118:66--71, 2011.

\bibitem{Rifkin}
J.~Rifkin.
\newblock {\em The third industrial revolution}.
\newblock Palgrave Macmillan, 2011.

\bibitem{worldfinancialreview}
J.~Rifkin.
\newblock The third industrial revolution: How the internet, green electricity
  and 3d printing are ushering in a sustainable era of distributed capitalism.
\newblock {\em World Financial Review}, 2012.

\bibitem{contemporary}
R.Q. R.Q.~Berry, G.Bull, C.Browning, D.D. Thomas, K.Starkweather, and J.H.
  Aylor.
\newblock Preliminary considerations regarding use of digital fabrication to
  incorporate engineering design principles in elementary mathematics
  education.
\newblock {\em Contemporary Issues in Technology and Teacher Education},
  10(2):167--172, 2010.

\bibitem{Saito}
K.~Saito.
\newblock Reading ancient greek mathematics.
\newblock In Eleanor Robson and Jacqueline Stedall, editors, {\em The Oxford
  Handbook of the history of Mathematics}. Oxford University Press, 2009.

\bibitem{Simms}
D.L. Simms.
\newblock The trail for {A}rchimedes's tomb.
\newblock {\em Journal of the Warburg and Courtauld Institutes}, 53:281--286,
  1990.

\bibitem{Slavkovsky}
E.~Slavkovsky.
\newblock {\em Feasability study for teaching geometry and other topics using
  three-dimensional printers}.
\newblock Harvard University, 2012.
\newblock A thesis in the field of mathematics for teaching for the degree of
  Master of Liberal Arts in Extension Studies.

\bibitem{Stein}
S.~Stein.
\newblock {\em Archimedes, what did he do besides Cry Eureka}.
\newblock The Mathematical Association of America, 1999.

\bibitem{Thomas}
I.~Thomas.
\newblock {\em Selections illustrating the history of Greek Mathematics}.
\newblock Harvard University Press,London, Heinemann, third edition edition,
  1957.

\bibitem{Trott}
M.~Trott.
\newblock {\em The Mathematica Guide book}.
\newblock Springer Verlag, 2004.

\bibitem{Turnbull}
H.W. Turnbull.
\newblock {\em The great mathematicians}.
\newblock Methuen and Co, Ltd, London, 1929.

\bibitem{Wohler}
T.~Wohlers.
\newblock 3d printing in education.
\newblock {\em Wohlerassociates}, 1995.
\newblock http://wohlersassociates.com/SepOct05TCT3dp.htm.

\end{thebibliography}

\end{document}